\documentclass[a4paper,fleqn]{cas-sc}

\usepackage[numbers]{natbib}

\usepackage{amsmath}
\numberwithin{equation}{section}
\usepackage{makecell}
\usepackage{diagbox}
\usepackage{subfigure}
\usepackage{color}
\usepackage{caption}
\graphicspath{{./figures/}}
\usepackage{graphicx}
\usepackage{appendix}
\usepackage{epstopdf}
\newcommand{\txrev}{\textcolor{black}}
\def\tsc#1{\csdef{#1}{\textsc{\lowercase{#1}}\xspace}}
\tsc{WGM}

\tsc{QE}
\begin{document}
\let\WriteBookmarks\relax
\def\floatpagepagefraction{1}
\def\textpagefraction{.001}

\shorttitle{Learning-based Multi-continuum Model for Numerical Homogenization}

\shortauthors{Fan Wang et al.}

\title [mode = title]{Learning-based Multi-continuum Model for Multiscale Flow Problems}

%\tnotemark[]

%\tnotetext[1]{tnote text}


\author[1]{Fan Wang}
\ead{wangfan525@stu.xjtu.edu.cn}

\author[1]{Yating Wang}
\cormark[1]
\ead{yatingwang@xjtu.edu.cn}

\author[2]{Wing Tat Leung}
\ead{wtleung27@cityu.edu.hk}

\author[1]{Zongben Xu}
\ead{zbxu@mail.xjtu.edu.cn}

\affiliation[1]{organization={School of Mathematics and Statistics, Xi'an Jiaotong University},city={Xi'an},country={China}}

\affiliation[2]{organization={Department of Mathematics, City University of Hong Kong},
            city={Hong Kong},
            country={Hong Kong Special Administrative Region}}

%\credit{Credit authorship details}

\cortext[1]{Corresponding author}



\begin{abstract}
%Numerical homogenization methods are classic ways of solving multiscale problems. They approximate a
 Multiscale problems can usually be approximated through numerical homogenization by an equation with some effective parameters that can capture the macroscopic behavior of the original system on the coarse grid to speed up the simulation. However, this approach usually assumes scale separation and that the heterogeneity of the solution can be approximated by the solution average in each coarse block. For complex multiscale problems, the computed single effective properties/continuum might be inadequate. In this paper, we propose a novel learning-based multi-continuum model to enrich the homogenized equation and improve the accuracy of the single continuum model for multiscale problems with some given data. Without loss of generalization, we consider a two-continuum case. The first flow equation keeps the information of the original homogenized equation with an additional interaction term. The second continuum is newly introduced, and the effective permeability in the second flow equation is determined by a neural network. The interaction term between the two continua aligns with that used in the Dual-porosity model but with a learnable coefficient determined by another neural network. The new model with neural network terms is then optimized using trusted data. We discuss both direct back-propagation and the adjoint method for the PDE-constraint optimization problem. Our proposed learning-based multi-continuum model can resolve multiple interacted media within each coarse grid block and describe the mass transfer among them, and it has been demonstrated to significantly improve the simulation results through numerical experiments involving both linear and nonlinear flow equations.
   % Numerical homogenization methods are classic ways of solving multiscale problems. They approximate a multiscale equation by a homogenized equation without multiscale properties, which can capture the macroscopic behavior of the original system. However, this approach can result in a relatively large error in the solution. In this paper, we propose a novel learning-based multi-continuum model to correct the homogenized equation and improve the accuracy of the numerical homogenization solution for multiscale problems. Specifically, the homogenized equation is reformulated as a multi-continuum model. The first continuum's flow equation is the homogenized equation obtained in numerical homogenization. The permeability of the second continuum's flow equation is determined by a neural network. The interaction term between the two continua aligns with those used in the Dual-porosity model but with a learnable coefficient determined by another neural network. These two networks are optimized using trusted data. Additionally, we have derived formulas for the gradient of the loss function with respect to the network parameters using the adjoint method. Our proposed learning-based multi-continuum model has been demonstrated to significantly improve the accuracy of numerical homogenization solutions through numerical experiments involving both linear and nonlinear flow equations.
\end{abstract}


%\item
%\item
%\item

% Keywords
\begin{keywords}
Multiscale problems \sep Numerical homogenization\sep Multi-continuum model \sep Deep learning 
\end{keywords}

\maketitle

% Main text
\section{Introduction}\label{intro}
Multiscale flow problems occur naturally across various physical and engineering applications. The subsurface formation of the flow problem usually involves discrete fractures and faults. Modeling and simulating processes within such fractured media presents difficulties due to the complexity of material properties, such as thermal diffusivity or hydraulic conductivity. The substantial differences between the fractures and the background pose challenges in numerical simulation since the high-contrast feature of the heterogeneous media introduces stiffness in the system and leads to a heavy computational burden for full fine-scale computation.

Multiscale model reduction methods, including local \cite{aarnes2008mixed1, aarnes2008mixed2, allaire2005multiscale,efendiev2011multiscale} and global approaches \cite{hinze2005proper, chinesta2011short, benner2015survey,chaturantabut2010nonlinear}, have been proposed to reduce computational costs in the numerical simulation of multiscale problems. The main idea is to construct reduced-order models to approximate the full fine-scale model to achieve efficient computation. Many local model reduction approaches, such as numerical homogenization \cite{BABUSKA197689,Nguetseng1989,Allaire1992,jikov2012homogenization,engquist2008asymptotic}, multiscale finite element methods \cite{efendiev2009multiscale, efendiev2013generalized}, variational multiscale methods \cite{hughes1998variational}, heterogeneous multiscale methods \cite{abdulle2012heterogeneous}, localized orthogonal decomposition \cite{maalqvist2014localization,henning2014localized,henning2013oversampling}, etc, have been successfully applied in many applications. %\cite{durlofsky1997nonuniform,yang2015upscaling}

Among these, numerical homogenization methods/flow-based upscaling methods are classic ways of solving multiscale problems. They calculate effective properties, such as permeability, in each coarse block that is comprised of an agglomeration of fine grid blocks, to incorporate fine-scale information into coarse-scale properties. Through the homogenization process, a complex multiscale system is approximated by a simple homogenized equation in which the heterogeneous microstructure is replaced by an equivalent homogeneous macrostructure. This homogenized equation can be numerically solved on a coarse mesh in a fast manner. 

However, the standard numerical homogenization is insufficient for multiple interacted media. A multiple continua approach is needed. The multi-continuum model \cite{barenblatt1960basic,abdassah1986triple,bai1993multiporosity,wu2006multiple} is proposed by considering multiple interconnected parallel systems that are distributed all over the domain in the flow field. The dual-porosity model \cite{barenblatt1960basic} is the first multi-continuum model in which the flow equations for the background (called the matrix) and the fracture are written separately with some interaction. In these approaches, several effective properties need to be computed locally in each coarse block, and some transfer terms are introduced to model the interaction. The transfer coefficients in the multi-continuum model are usually obtained based on simplified physical assumptions.

For complicated physical processes, the above-mentioned reduced order models might be inadequate due to ignoring microscopic information or physical simplifying assumptions. Thus resulting in non-negligible errors between the simulation results obtained from the reduced order models and the real observation data. In this paper, we will consider improving a given inaccurate reduced order model with trusted data by a deep learning approach. 

Recently, leveraging the deep learning approach to assist scientific computing has been an emerging trend. Many works have been done to apply neural networks to approximate the solutions of partial differential equations (PDEs), such as PINN~\cite{raissi2019physics}, deep Ritz method~\cite{yu2018deep}, deep Galerkin method~\cite{sirignano2018dgm}, etc. For parameterized PDEs, operator learning approaches such as DeepOnet \cite{lu2021learning} and Fourier Neural Operator \cite{li2020fourier} have been developed to learn the mapping from the parameter space to the solutions space. Other deep learning approaches for PDE-related tasks include developing surrogate models motivated by classical numerical solvers. For example, in \cite{Fan2018MNNH, PCDL_nz}, the authors construct multiscale neural networks based on hierarchical multigrid solvers and encoder-decoder neural networks for solutions of heterogeneous elliptic PDEs. There are also other deep learning approaches to handle multiscale problems. For instance, some data-driven methods are proposed to obtain the coarse grid effective properties based on the nonlocal multi-continuum upscaling method \cite{wang2018NLMC_DNN}, and the coefficients in global model reduction techniques such as the proper orthogonal decomposition (POD) projection \cite{pod_dl_SW}. Furthermore, a neural homogenization-based PINN (NH-PINN) \cite{LEUNG2022111539} was investigated to improve the PINN accuracy for solving multiscale problems with the help of homogenization. The idea is to use PINN to solve the cell problems and the homogenized coefficients can be evaluated to form the homogenized equation. Then, PINN is utilized again to solve the global homogenized equation. However, in the above methods, the physical model is assumed to be accurate, and deep learning methods are used to learn the solution of the underlying physical model.

Another area of research involves representing unknown or omitted physical processes by deep learning methods, in order to correct the original PDE to make them more accurate and easier to solve. These unrepresented or unknown physics are typically challenging to express due to limited cognition or computational cost. The deep learning PDE augmentation method (DPM) \cite{SIRIGNANO2020109811} leverages known physics (expressed in a PDE) to learn closures for unknown or unrepresented physics. For the large-eddy simulation of turbulence, DPM indicates that the neural network can acquire a more effective sub-filter-scale closure than conventional LES models, which enhances numerical accuracy. The approaches in \cite{MILLER2022111541, XIAO2023112317} employ neural networks to approximate the collision operator of the Boltzmann equation. Because of the extremely high dimensionality and nonlinearity, this collision operator poses a notorious challenge for theoretical analysis and numerical simulation. Compared with the classical simplified model, the neural network-based model of a collision operator not only agrees with the direct simulation but also maintains the essential structural properties of a many-particle system. The study in \cite{SIRIGNANO2023112016} presents a theoretical analysis of the convergence of optimizing neural network terms of linear elliptic PDEs by the adjoint method.

In this paper, we proposed a novel learning-based multi-continuum model, to enrich the homogenized equation and address the inaccuracy of PDE solutions caused by numerical homogenization for multiscale problems. Without loss of generalization, we consider a two-continuum case. We remodel the homogenized single continuum equation as a multi-continuum model, where the homogenized equation serves as the flow equation for the first continuum with an additional interaction term. For the newly introduced second continuum, the effective permeability of its flow equation is determined by a neural network. The interaction term of these two continua aligns with that used in the Dual-porosity model but with a learnable coefficient determined by another neural network. By optimizing networks through the trusted data, our learning-based multi-continuum model improves the accuracy of solutions for the original homogenized equation. Additionally, we discuss both direct back-propagation and the adjoint method for the PDE-constraint optimization problem. The main contributions of this paper are summarized as follows.

% In this paper, we proposed a novel learning-based multi-continuum model, to address the inaccuracy of PDE solutions caused by numerical homogenization for multiscale problems. We remodel the homogenized single continuum equation as a multi-continuum model, where the homogenized equation serves as the flow equation for the first continuum with an additional interaction term. For the second continuum's flow equation, its effective permeability is determined by a neural network. The interaction term of these two continua aligns with that used in the Dual-porosity model but with a learnable coefficient determined by another neural network. By optimizing networks through the trusted data, our learning-based multi-continuum model improves the accuracy of solutions for the original homogenized equation. Additionally, we derive formulas for the gradient of the loss function with respect to the network parameters using the adjoint method. The main contributions of this paper are summarized as follows.

\begin{itemize}
    \item We propose a new learning-based multi-continuum model to correct the homogenized flow equation with given data and improve the accuracy of the solution for multiscale problems.
    \item We consider both direct back-propagation and the adjoint method to compute the gradients of the loss function with respect to the network parameters to optimize our proposed learning-based multi-continuum model.
    \item We conducted numerical experiments including linear and nonlinear flow equations. Our proposed model can substantially improve the accuracy of the solution of the reduced order model obtained from standard numerical homogenization.
\end{itemize}

The rest of the paper is organized as follows. In Section \ref{back}, we review the basics of the multiscale flow problem, the reduced-order model with numerical homogenization, and the multi-continuum model. Next, we will present our proposed learning-based multi-continuum model in Section \ref{Method}. We conduct several numerical experiments in Section \ref{experiments}. Conclusions are given in Section \ref{conclusion}.

\section{Background}\label{back}
\subsection{Multiscale flow problem}
Consider the parabolic equation
\begin{equation}\label{multi-scale}
\begin{aligned}
    \frac{\partial u(\mathbf{x},t)}{\partial t} - \text{div}(\kappa(\mathbf{x},u(\mathbf{x},t))\nabla u(\mathbf{x},t)) &= f(\mathbf{x},t), &\text{in} \ \Omega \times [0,T],\\
    u(\mathbf{x}, t) &= 0, &\text{in} \  \Omega \times \{0\},\\
    u(\mathbf{x}, t) &= 0, &\quad \text{on} \  \partial\Omega\times [0,T],
\end{aligned}
\end{equation}
where $u(\mathbf{x},t)\in \mathbb{R}$ is the solution, $\mathbf{x}\in \Omega \subset \mathbb{R}^d$, $T>0$ is the final time.%, and $\kappa(\mathbf{x},u(\mathbf{x},t))\in \mathbb{R}^{d\times d}$ is the value of the conductivity/permeability. 
\txrev{The nonlinear diffusion coefficient $\kappa(\mathbf{x},u(\mathbf{x},t))$ represents the conductivity/permeability of the underlying field}, and we assume $\kappa(\mathbf{x},u(\mathbf{x},t)) = \kappa_0(\mathbf{x}) \kappa_u(u(\mathbf{x},t))$. In particular, assume $\kappa_u(u(\mathbf{x},t)) > 0$ and is bounded and Lipchitz continuous, and $\kappa_0 \in L^{\infty}(\Omega)$ is the value of the conductivity/permeability of a fractured/channelized media. Since the material properties within fractures can be very different from the background properties and the complex fracture configuration, there are multiple scales and high contrast of $\kappa_0$. The multiscale and nonlinear feature poses the challenges in the numerical simulation.% due to the high-contrast $\kappa$ of the heterogeneous media.

\subsection{Reduced order model with numerical homogenization -- single continuum}\label{ROM}
Numerical homogenization is a classic model reduction approach to finding numerical solutions of heterogeneous multiscale problems. It approximates the original multiscale problem with a simpler homogenized one. The solution of this homogenized equation can capture the macroscopic behavior of the original multiscale problem. The parameters of the homogenized equation (such as permeability) are normally homogeneous or smooth, such that the homogenized equation can be solved numerically in a coarse mesh.

To illustrate the idea, we here give the numerical homogenization of the parabolic equation with the coefficients that has two scales $\mathbf{x}$ and $\mathbf{y} = \mathbf{x}/\epsilon$, 
\txrev{
\begin{equation}
    \frac{\partial u_\epsilon}{\partial t} - \frac{\partial}{\partial x_i}(\kappa_{ij}(\mathbf{x}, \mathbf{x}/\epsilon) \frac{\partial}{\partial x_j} )u_\epsilon = f,\quad (\mathbf{x},t) \in \Omega \times [0,T],\\
\end{equation}}
with $u_\epsilon(\mathbf{x}, t) = 0$ on $\partial \Omega$ or $t=0$. \txrev{Here we use the Einstein notation}. We seek for $u_\epsilon(\mathbf{x}, t)$ in the asymptotic expansion
\begin{equation}
    u_\epsilon(\mathbf{x},t)=u_0(\mathbf{x},\mathbf{x}/\epsilon,t)+\epsilon u_1(\mathbf{x},\mathbf{x}/\epsilon,t) + \epsilon^2 u_2(\mathbf{x},\mathbf{x}/\epsilon,t) + \mathcal{O}(\epsilon^3),
\end{equation}
\txrev{where $u_j(\mathbf{x},\mathbf{x}/\epsilon,t)$, for $j=0,1,2,\ldots,$ are periodic in $\mathbf{y}=\mathbf{x}/\epsilon$}. Denoted by $A^\epsilon$ the second-order elliptic operator
\begin{equation}
    A^\epsilon = - \frac{\partial}{\partial x_i}(\kappa_{ij}(\mathbf{x}, \mathbf{x}/\epsilon) \frac{\partial}{\partial x_j} ).
\end{equation}
With this notation, it can be checked that
\begin{equation}
    A^\epsilon = \epsilon^{-2}A_1 + \epsilon^{-1}A_2 + \epsilon^0 A_3,
\end{equation}
where
\begin{equation}
\begin{aligned}
& A_1 = - \frac{\partial}{\partial y_i}(\kappa_{ij}(\mathbf{x}, \mathbf{y}) \frac{\partial}{\partial y_j} ),\\
& A_2 = - \frac{\partial}{\partial y_i}(\kappa_{ij}(\mathbf{x}, \mathbf{y}) \frac{\partial}{\partial x_j} ) - \frac{\partial}{\partial x_i}(\kappa_{ij}(\mathbf{x}, \mathbf{y}) \frac{\partial}{\partial y_j}),\\
& A_3 = - \frac{\partial}{\partial x_i}(\kappa_{ij}(\mathbf{x}, \mathbf{y}) \frac{\partial}{\partial x_j} ).
\end{aligned}
\end{equation}
We hence have $A^\epsilon u_\epsilon=f$. By equating the terms of different powers of $\epsilon$, we have
\begin{equation}\label{eq:nh_para}
\begin{aligned}
    & A_1u_0 = 0,\\
    & A_1u_1 + A_2u_0 = 0,\\
    & A_1u_2 + A_2u_1 + A_3u_0 -f = -\frac{\partial u_0}{\partial t}.
\end{aligned}
\end{equation}
Substitute $A_1$ into the first equation of \eqref{eq:nh_para}, we can obtain 
\begin{equation}
- \frac{\partial}{\partial y_i}(\kappa_{ij}(\mathbf{x}, \mathbf{y}) \frac{\partial}{\partial y_j} ) u_0(\mathbf{x},\mathbf{y},t) = 0
\end{equation}
where $u_0$ is periodic in $\mathbf{y}$. The theory of second-order elliptic PDEs implies that $u_0(\mathbf{x},\mathbf{y},t)$ is independent of $\mathbf{y}$, i.e., $u_0(\mathbf{x},\mathbf{y},t) = u_0(\mathbf{x},t)$. This leads to 
\begin{equation}
- \frac{\partial}{\partial y_i}(\kappa_{ij}(\mathbf{x}, \mathbf{y}) \frac{\partial}{\partial y_j} )u_1 = (\frac{\partial}{\partial y_i}\kappa_{ij}(\mathbf{x}, \mathbf{y})) \frac{\partial u_0(\mathbf{x},t)}{\partial x_j}.
\end{equation}
Define $\mathcal{N}_j = \mathcal{N}_j (\mathbf{x},\mathbf{y}) $ as the solution to the following \textit{cell problem}:
\begin{equation}
- \frac{\partial}{\partial y_i}(\kappa_{ij}(\mathbf{x}, \mathbf{y}) \frac{\partial}{\partial y_j} )\mathcal{N}_j = \frac{\partial}{\partial y_i}\kappa_{ij}(\mathbf{x}, \mathbf{y}),
\end{equation}
\txrev{where $\mathcal{N}_j$ is periodic in $\mathbf{y}$}. The general solution for $u_1$ is then given by
\txrev{
\begin{equation}
u_1(\mathbf{x},\mathbf{y},t) = \mathcal{N}_j(\mathbf{x},\mathbf{y})\frac{\partial u_0(\mathbf{x},t)}{\partial x_j}.
\end{equation}}
Finally, we note that the eqaution for $u_2$ is given by
\txrev{
\begin{equation}\label{1}
\frac{\partial}{\partial y_i}(\kappa_{ij}(\mathbf{x}, \mathbf{y}) \frac{\partial}{\partial y_j} ) u_2 = A_2u_1+A_3u_0-f +\frac{\partial u_0}{\partial t}.
\end{equation}}
The solvability condition implies that the right-hand side of (\ref{1}) must have mean zero in $\mathbf{y}$ over one periodic cell $Y=[0,1]\times[0,1]$, i.e.,
\begin{equation}
\int_Y (A_2u_1+A_3u_0-f +\frac{\partial u_0}{\partial t}) d\mathbf{y}=0.
\end{equation}
Next, a different trick is used. In particular, the equation at $\epsilon^0$ is averaged over the period. We note that
 \txrev{
\begin{equation}
\langle \frac{\partial}{\partial y_i}F(\mathbf{x},\mathbf{y},t) \rangle = \int_Y \frac{\partial}{\partial y_i} F(\mathbf{x},\mathbf{y},t) d\mathbf{y}=0,
\end{equation}}
\txrev{for any $F(\mathbf{x},\mathbf{y}, t)$ which is periodic with respect to $\mathbf{y}$}. This can be easily verified using the divergence theorem. \txrev{Then, the averages over the terms which start with $\frac{\partial}{\partial y_i}$ will disappear and we get}
\txrev{
\begin{equation}
\frac{\partial u_0}{\partial t}-\frac{\partial}{\partial x_i} (\langle \kappa_{ij}(\mathbf{x},\mathbf{y})\rangle \frac{\partial}{\partial x_j}u_0) - \frac{\partial}{\partial x_i} \langle \kappa_{ij}(\mathbf{x}, \mathbf{y})\frac{\partial}{\partial y_j}u_1\rangle  = f. 
\end{equation}}

Substituting the expression for $u_1$ into this equation, we obtain
\begin{equation}
    \frac{\partial u_0}{\partial t} - \frac{\partial}{\partial x_i}(\kappa^\star_{ij} \frac{\partial}{\partial x_j} )u_0 = f, 
\end{equation}
where 
\txrev{ 
\begin{equation}
    \kappa^\star_{ij} = \frac{1}{|Y|} \int_Y (\kappa_{ij}+\kappa_{ik}\frac{\partial \mathcal{N}_j}{\partial y_k}) d\mathbf{y}.
\end{equation}
}
The homogenized equation of the multiscale flow problem (\ref{multi-scale}) is then given by
\txrev{
\begin{equation}\label{homogenization}
      \frac{\partial u(\mathbf{x},t)}{\partial t} - \text{div}(\kappa^\star \nabla u(\mathbf{x},t)) = f(\mathbf{x},t) \quad \text{in} \quad \Omega \times [0,T],
\end{equation}}
where $\kappa^\star$ is the homogenized value of the  permeability. \txrev{The homogenized solution is a kind of local average of the solution of the original multiscale equation. This homogenized solution approximates the original solution and can capture the macroscopic behavior of the original multiscale problem. We build on the homogenized equation as a good initialization and further correct it to obtain a more accurate model for characterizing the complex dynamics of multiscale problems.}

\subsection{Multi-continuum model.}
The numerical homogenization technique is widely used for single continuum upscaling. However, in a more general setup, a multi-continuum approach is needed. The multi-continuum model divides a material or medium into multiple regions with different physical properties between each region. These regions can have different structures, compositions, porosities, temperatures, and so on. The model describes the behavior of the whole system by defining the specific physical properties of each region (for example, fluid flow, heat transfer, mass transport, etc.) and the interactions between them. By incorporating various physical processes and considering the mutual coupling between regions, multi-continuum models can more accurately capture the macroscopic behavior of a complicated material or medium.

A coupled system of dual-continuum equations for the multscale flow equation has the form
\begin{equation}\label{multi-continuum}
  \begin{aligned}
      \frac{\partial u_1(\mathbf{x},t)}{\partial t} - \text{div}(\kappa_1(\mathbf{x},u_1(\mathbf{x},t))\nabla u_1(\mathbf{x},t)) + \sigma_{12}(\mathbf{x},u_1(\mathbf{x},t), u_2(\mathbf{x},t))(u_1(\mathbf{x},t)-u_2(\mathbf{x},t)) = f_1(\mathbf{x},t), \\
       \frac{\partial u_2(\mathbf{x},t)}{\partial t} - \text{div}(\kappa_2(\mathbf{x},u_2(\mathbf{x},t))\nabla u_2(\mathbf{x},t)) + \sigma_{21}(\mathbf{x},u_2(\mathbf{x},t), u_1(\mathbf{x},t))(u_2(\mathbf{x},t)-u_1(\mathbf{x},t)) = f_2(\mathbf{x},t),
  \end{aligned}
\end{equation}
where $(\mathbf{x},t)\in\Omega\times[0,T]$, $u_i$ denotes the solution for the $i$-th continuum, $\kappa_i$ denotes the vale of permeability, $f_i$ represents the source function for the $i$th continuum, $i=1,2$. Terms $\sigma_{12}(\mathbf{x},u_1, u_2)(u_1-u_2)$ and $\sigma_{21}(\mathbf{x},u_2, u_1)(u_2-u_1)$ describe mass transfer of the liquid which flows from one continuum into another continuum per unit of media volume as well as per unit of time. However, it is usually difficult to capture the transient fluid transfer between different continua, and the transfer coefficients $\sigma_{12}$ and $\sigma_{21}$ are typically determined empirically.

\section{Method}\label{Method}
\begin{figure}[t]
        \centering
        \includegraphics[width=1\textwidth]{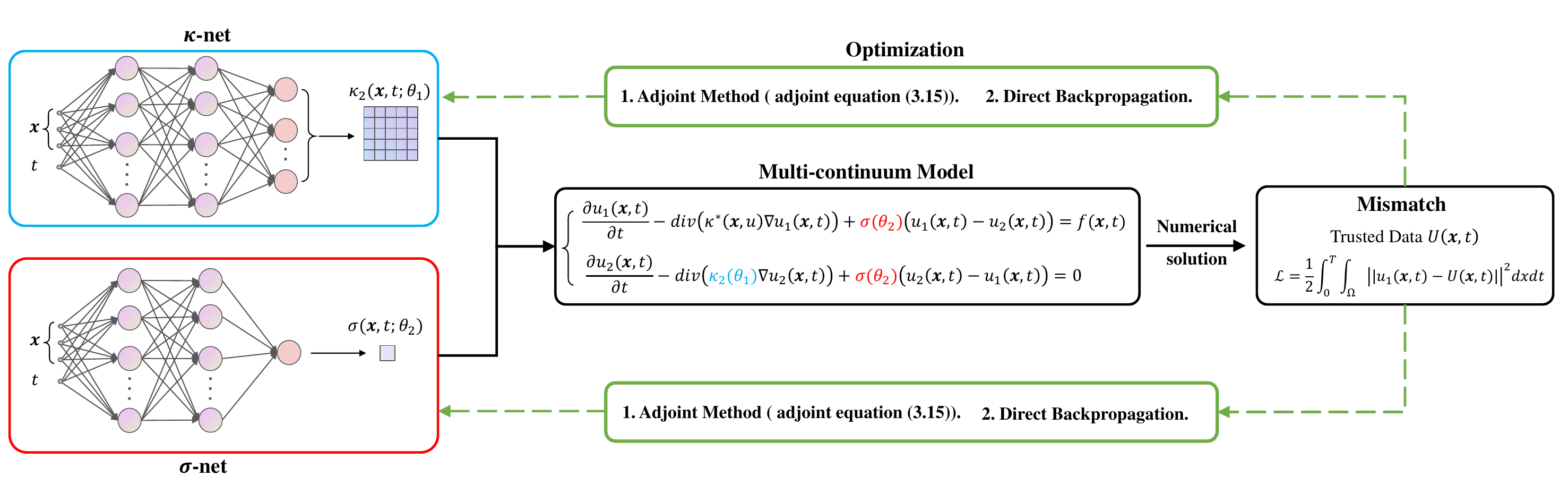}
        \caption{The schematic illustration of our proposed learning-based multi-continuum model. The black color indicates the forward process, and the green color indicates the back-propagation to optimize network parameters.}
        \label{method}
\end{figure}

In this section, we present our proposed approach to model correction based on the deep learning method. We propose to rectify the inaccuracy of numerical homogenization of the highly heterogeneous multiscale flow equation, especially the ones that possess multi-continuum properties in the underlying media, using a learning-based multi-continuum model. Our model involves learning the permeability in the new continuum and transfer coefficient from the trusted data. \textbf{Figure} \ref{method} illustrates our methodology. The permeability and transfer coefficient are parameterized as neural networks. The network parameters can be optimized by aligning the numerical solutions of the corrected model with trusted data. To optimize the networks, we have selected two methods to calculate the gradients of the network: direct back-propagation (BP) with the automatic differentiation of Tensorflow/Pytorch and the adjoint method. The adjoint method computes the loss function's derivatives relative to the network parameters by solving the adjoint equation, thus eliminating the need for abundant storage in the forward solver. However, computing the solutions to the adjoint equations also introduces additional computational costs. In the following subsections, we introduce our learning-based multi-continuum model, as well as the forward PDE solvers and optimization methods for linear and non-linear problems, respectively.

\subsection{Learning-based multi-continuum model}\label{sec3.1}
Numerical homogenization captures the macroscopic behavior of the original complicated multiscale flow equation and generates solutions with low precision. To increase the accuracy of solutions, we propose to use a multi-continuum model to enrich the homogenized equation (\ref{homogenization}), which is formulated as
\begin{equation}\label{Correction}
    \begin{aligned}
    &\frac{\partial u_1(\mathbf{x},t)}{\partial t} - \text{div}(\kappa_1^\star(\mathbf{x},u_1)\nabla u_1(\mathbf{x},t)) + \sigma(\mathbf{x},u_1,u_2)(u_1(\mathbf{x},t)-u_2(\mathbf{x},t)) = f(\mathbf{x},t),\  & \text{in} \ & \Omega \times [0,T],\\
    &\frac{\partial u_2(\mathbf{x},t)}{\partial t} - \text{div}(\kappa_2(\mathbf{x},u_2)\nabla u_2(\mathbf{x},t)) + \sigma(\mathbf{x},u_1,u_2)(u_2(\mathbf{x},t)-u_1(\mathbf{x},t)) = 0,\  & \text{in} &\  \Omega \times [0,T],
    \end{aligned}
\end{equation}
where $u_1$ corresponds to the solution of the homogenized flow equation, $\kappa_1^\star$ denotes the homogenized permeability, $u_2$ is the solution of the flow equation in the introduced second continuum, $\kappa_2(\mathbf{x}, u_2)$ denotes the permeability of the second continuum, and $\sigma(\mathbf{x},u_1,u_2)$ is the transfer coefficient. The source term for the second continuum's flow equation is set to 0. The initial and boundary conditions of this multi-continuum model are given by
$$ u_1(\mathbf{x}, t=0) = u_2(\mathbf{x}, t=0) = 0, \quad \text{in} \  \Omega, \qquad
    u_1(\mathbf{x}, t) = u_2(\mathbf{x}, t) = 0, \quad \text{on} \  \partial\Omega\times [0,T]. $$

In model (\ref{Correction}), both the permeability field $\kappa_2$ and the transfer coefficient $\sigma$ are unknown and need to be determined. For the newly introduced continuum, we do not know its properties. Thus we propose to leverage the deep learning method to determine the values of the permeability field $\kappa_2$ and the transfer coefficient $\sigma$. The schematic illustration of our method is shown in \textbf{Figure} \ref{method}. Specifically, the learning-based multi-continuum model is defined as
\begin{equation}\label{Learning-Correction}
    \begin{aligned}
    &\frac{\partial u_1(\mathbf{x},t)}{\partial t} - \text{div}(\kappa_1^\star(\mathbf{x},u_1)\nabla u_1(\mathbf{x},t)) + \sigma(\mathbf{x},t;\theta_2)(u_1(\mathbf{x},t)-u_2(\mathbf{x},t)) = f(\mathbf{x},t),\  & \text{in} &\  \Omega \times [0,T],\\
    &\frac{\partial u_2(\mathbf{x},t)}{\partial t} - \text{div}(\kappa_2(\mathbf{x},t;\theta_1)\nabla u_2(\mathbf{x},t)) + \sigma(\mathbf{x},t;\theta_2)(u_2(\mathbf{x},t)-u_1(\mathbf{x},t))= 0,\  & \text{in} &\  \Omega \times [0,T],
    \end{aligned}
\end{equation}
where $\kappa_2(\mathbf{x},t;\theta_1)$ and $\sigma(\mathbf{x},t;\theta_2)$ are parameterized as neural networks, named $\kappa$-net and $\sigma$-net respectively, $\theta_1$ and $\theta_2$ are network parameters. The $\kappa$-net $\kappa_2(\mathbf{x}, t;\theta_1):\mathbb{R}^{d+1}\rightarrow \mathbb{R}^{d\times d}$ learns elements of the permeability matrix of the second continuum at each coordinate $\mathbf{x}$ and time $t$. It corrects the solution $u_1$ by controlling the dynamics of $u_2$ in the second newly introduced continuum. \txrev{The $\sigma$-net $\sigma(\mathbf{x}, t; \theta_2):\mathbb{R}^{d+1}\rightarrow \mathbb{R}$ learns a transfer coefficient to determine the impact of $u_2$ on the solution $u_1$.}

In this way, we obtain a learning-based multi-continuum model that corrects the original homogenized equation. Its numerical solution can be obtained by classical methods (such as finite difference method (FDM) and finite element method (FEM)) and deep learning methods (PINN, DeepRitz, etc.). The network parameters $\theta_1$ and $\theta_2$ are optimized by minimizing the mismatch between the numerical solution and the trusted data, i.e. minimizing the loss function
% \textbf{Optimization problem}
\begin{equation}\label{loss}
    \mathcal{L}=\frac{1}{2}\int_0^T \int_\Omega \| u_1(\mathbf{x}, t) - U(\mathbf{x}, t) \|^2 d \mathbf{x}dt,
\end{equation}
where $U(\mathbf{x},t)$ denotes the trusted data of the multiscale flow equation (\ref{multi-scale}) assumed to be available at certain coordinates and times. By minimizing the loss function (\ref{loss}), the permeability field $\kappa_2$ and transfer coefficient $\sigma$ can be determined and we obtain a new multi-continuum model that more accurately models the dynamics of the flow in the heterogeneous media.

%\begin{equation}\label{loss}
%    L=\sum\limits_{m=1}^M \sum\limits_{n=1}^{N_t} \| u(x_m, t_n) - U(x_m, t_n) \|_2^2,
%\end{equation}
%where $\{u(x_m, t_n)\}_{m=1,n=1}^{M,N_t}$ is the numerical solutions of equation (\ref{Correction}) and $U(t,x)$ are trusted data assumed to be available at certain times $\{t_n\}_{n=1}^{N_t}$ and coordinate $\{x_m\}_{m=1}^{M}$. By minimizing the above loss function, a new two-continuum model can be learned to corrected the original homogenized equation, and it can generates the solution $u$ as close as possible to the trust data $U$.

\subsection{Forward solver}\label{forward}
The calculation of the loss function (\ref{loss}) requires the numerical solutions of equation (\ref{Learning-Correction}). In this section, we present our numerical solver for linear and nonlinear equations. FEM is the basic solver used in our method. In addition, we propose to use the deep learning method, PINN, to solve the nonlinear equation to simplify the arduous iterations. Details are as follows.

%The optimization process of the network parameters involves solving equation (\ref{Correction}) as discussed in section \ref{sec3.1}. The finite element method is the most basic forward solver. While for nonlinear flow equation, using Newton's method to solve the equation iteratively leads to a computationally cumbersome optimization process of network parameters. Thus we propose to use deep learning based method as forward solver, such as PINN\cite{}, DeepRitz\cite{}, etc. Next we will introduce the forward solvers for linear and nonlinear flow equation respectively.

\subsubsection{Linear flow equation}
For the linear flow equation, we assume the values of permeability only dependent on the coordinate $\mathbf{x}$, i.e.,
\begin{equation}\label{linear}
    \frac{\partial u(\mathbf{x},t)}{\partial t} - \text{div}(\kappa(\mathbf{x})\nabla u(\mathbf{x},t)) = f(\mathbf{x},t), \quad \text{in} \quad \Omega \times [0,T].
\end{equation}
The corresponding learning-based two-continuum correction model is reformulated as
\begin{equation}\label{Correction-linear}
    \begin{aligned}
    & \frac{\partial u_1(\mathbf{x}, t)}{\partial t} - \text{div}(\kappa_1^\star(\mathbf{x})\nabla u_1(\mathbf{x},t)) + \sigma(\mathbf{x};\theta_2)(u_1(\mathbf{x},t)-u_2(\mathbf{x},t)) = f(\mathbf{x},t),\  & \text{in} &\  \Omega \times [0,T],\\
    & \frac{\partial u_2(\mathbf{x,t})}{\partial t} - \text{div}(\kappa_2(\mathbf{x};\theta_1)\nabla u_2(\mathbf{x},t)) + \sigma(\mathbf{x};\theta_2)(u_2(\mathbf{x},t)-u_1(\mathbf{x},t))= 0,\  & \text{in} &\  \Omega \times [0,T],
    \end{aligned}
\end{equation}
where the permeability $\kappa_2$ and transfer coefficient $\sigma$ are only dependent on the coordinate $\mathbf{x}$. We solve equation (\ref{Correction-linear}) by FEM. The weak form of equation (\ref{Correction-linear}) is to seek $u_1(t,\cdot), u_2(t,\cdot) \in V = H_0^1(\Omega)$ such that for $\forall v_1,v_2 \in V$,
\txrev{
\begin{equation}\label{weak}
\begin{aligned}
  & (\frac{\partial u_1}{\partial t}, v_1) + a_1(u_1, v_1) + (\sigma(u_1-u_2), v_1)= (f, v_1), \  &t\in(0, T],\\
  & (\frac{\partial u_2}{\partial t}, v_2) + a_2(u_2, v_2) + (\sigma(u_2-u_1), v_2)= 0, \  &t\in(0, T],
\end{aligned}
\end{equation}
where $a_1(u_1,v_1)=\int_\Omega \kappa_1^\star \nabla u_1 \cdot \nabla v_1 d\mathbf{x}$ and $a_2(u_2,v_2)=\int_\Omega \kappa_2 \nabla u_2\cdot \nabla v_2 d\mathbf{x}$.} Now consider a spatial partition $\Gamma_h$ of the computational domain $\Omega$, $h$ denotes the mesh size. We choose the linear basis function $\psi_i$ and form the space $V_h= \text{span}\{\psi_i, 1\leq i \leq L \}$. Then the discretization of (\ref{weak}) reads as
\txrev{
\begin{equation}\label{FE-linear}
\begin{aligned}
 & M\frac{\partial u_1}{\partial t}+ A_1(\kappa_1^\star)u_1 +\sigma(\mathbf{x}, \theta_2)M(u_1-u_2) = F,\\
 & M\frac{\partial u_2}{\partial t}+ A_2(\kappa_2(\mathbf{x}, \theta_1))u_2 +\sigma(\mathbf{x}, \theta_2)M(u_2-u_1) = 0,
\end{aligned}
\end{equation}}
where $M$ is the mass matrix whose elements are $\int_\Omega \psi_i \psi_j d\mathbf{x}$, $A_1$ and $A_2$ are the stiffness matrices with elements $\int_\Omega \kappa_1^\star \nabla \psi_i \cdot\nabla\psi_j d\mathbf{x}$ and $\int_\Omega \kappa_2 \nabla \psi_i\cdot \nabla\psi_j d\mathbf{x}$ respectively, and the elements of $F$ are $\int_\Omega f\psi_j d\mathbf{x}$. We use implicit backward Euler scheme in time
\begin{equation}
\begin{aligned}
 & M\frac{u_1^{k+1} - u_1^{k}}{\tau}+ A_1(\kappa^\star)u_1^{k+1} +\sigma(\mathbf{x}, \theta_2)M(u_1^{k+1}-u_2^{k+1}) = F^{k+1},\\
 & M\frac{u_2^{k+1} - u_2^{k}}{\tau}+ A_2(\kappa_2(\mathbf{x}, \theta_1))u_2^{k+1} +\sigma(\mathbf{x}, \theta_2)M(u_2^{k+1}-u_1^{k+1}) = 0,
\end{aligned}
\end{equation}
where $u^{k} = u(t_k)$ and $t_k=k\tau, k=0,1,\cdots,N_t-1$, $N_t=T/\tau$ is the number of time steps. For each time step $k$, the numerical solution of $u_1$ can be obtained by solving the above linear system.
%\textbf{forward solver: finite element method or PINN}

\subsubsection{Nonlinear equation}\label{nonlinear}
For the nonlinear multiscale flow equation (\ref{multi-scale}), we introduce two different methods, the classical FEM with iteration approach and the deep learning method PINN, to solve the learning-based multi-continuum model (\ref{Learning-Correction}).

\textbf{Finite element method}. The discretization of equation (\ref{Learning-Correction}) in FEM with implicit backward Euler scheme in time is formulated as
\begin{equation}\label{FE-nonlinear}
\begin{aligned}
 & M\frac{u_1^{k+1} - u_1^{k}}{\tau}+ A_1^{k+1}(\kappa_1^\star)u_1^{k+1} +\sigma^{k+1}(\mathbf{x},t, \theta_2)M(u_1^{k+1}-u_2^{k+1}) = F^{k+1},\\
 & M\frac{u_2^{k+1} - u_2^{k}}{\tau}+ A_2^{k+1}(\kappa_2(\mathbf{x},t, \theta_1))u_2^{k+1} +\sigma^{k+1}(\mathbf{x},t, \theta_2)M(u_2^{k+1}-u_1^{k+1}) = 0,
\end{aligned}
\end{equation}
where $A_1^{k+1}$ is dependent on $u^{k+1}_1$ and are different for each $k=0,2,\cdots,N_t-1$. The equation of $u_2$ is linear and $A_2^{k+1}$ we learned only determined by coordinate $\mathbf{x}$ and time $t$. In this paper, we use the Newton-Raphson method to solve the above nonlinear equation. This procedure employs a first order approximation of (\ref{FE-nonlinear}) in the neighbourhood of an approximate solution $u_1^{k+1;n}$ resulting in a linear system. For an approximate solution $u_1^{k+1;n}, n=1,2,...N$, we calculate the solution $u_2^{k+1;n}$ by
\begin{equation}
\begin{aligned}
 M\frac{u_2^{k+1;n} - u_2^{k;n}}{\tau}+ A_2^{k+1}(\kappa_2(\mathbf{x},t,\theta_1))u_2^{k+1;n} +\sigma^{k+1}(\mathbf{x},t,\theta_2)M(u_2^{k+1;n}-u_1^{k+1;n}) = 0,
\end{aligned}
\end{equation}
Then the approximate solution $u_1^{k+1;n}$ is updated according to
\begin{equation}
\begin{aligned}
 u_1^{k+1;n+1} = u_1^{k+1;n} + \delta u_1,
\end{aligned}
\end{equation}
\txrev{where $\delta u_1$ can be computed by} 
\begin{equation}
\begin{aligned}
 & \big(\frac{1}{\tau}M + \partial_{u_1^{k+1;n}}A_1^{k+1;n}(\kappa_1^\star)u_1^{k+1;n} + A_1^{k+1;n}(\kappa_1^\star) +\sigma^{k+1}(\mathbf{x},t,\theta_2)M \big)\delta u_1= -R_1^{k+1;n},
\end{aligned}
\end{equation}
where $A_1^{k+1;n}$ is the stiff matrix at $u_1^{k+1;n}$, and $R_1^{k+1;n}$ is the residuals of (\ref{FE-nonlinear}) at $u_1^{k+1;n}, u_2^{k+1;n}$. These so-called iterations are repeated until the residual of equation (\ref{FE-nonlinear}) is satisfied with sufficient precision.

Although the nonlinear multi-continuum model (\ref{Learning-Correction}) can be solved by utilizing the FEM and Newton-Raphson method, the latter involves numerous iterations which in turn result in high computational cost. Additionally, the backpropagation and adjoint method used to optimize the network parameters, as introduced in Section \ref{optimize}, also incur significant computation costs, due to multiple iterations in the Newton-Raphson method. Therefore, we have opted to employ the deep learning method, PINN, to solve the nonlinear equation (\ref{Learning-Correction}).

\textbf{Physics-informed neural network (PINN).} PINN is adopted as our forward solver for the nonlinear multiscale flow equation. It directly maps from the coordinate $\mathbf{x}$ and time $t$ to the desired solution. We denote the PINN solver as $\text{PINN}(\theta_3)$ where $\theta_3$ is the network parameter. By this way, the parameters $\theta_1, \theta_2$ of $\kappa$-net, $\sigma$-net , and parameters $\theta_3$ of PINN need to be optimized simultaneously. The loss function is reformulated as
\txrev{
\begin{equation}\label{pinn}
\begin{aligned}
 \mathcal{L}(\theta_1,\theta_2,\theta_3)&= \frac{1}{2}\int_0^T \int_\Omega \| u_1(\mathbf{x}, t) - U(\mathbf{x}, t) \|^2 d\mathbf{x}dt\\
            &+\int_0^T \int_\Omega \| \frac{\partial u_1}{\partial t} - \text{div}(\kappa_1^\star\nabla u_1) + \sigma(\mathbf{x},t;\theta_2)(u_1-u_2) - f\|^2d\mathbf{x}dt \\
            &+\int_0^T \int_\Omega \| \frac{\partial u_2}{\partial t} - \text{div}(\kappa_2(\mathbf{x},t;\theta_1)\nabla u_2) + \sigma(\mathbf{x},t;\theta_2)(u_2-u_1)\|^2d\mathbf{x}dt, \\
\end{aligned}
\end{equation}}
where $u_1,u_2 = \text{PINN}(\mathbf{x},t;\theta_3)$ are solutions obtained by the PINN solver. The PINN solver avoids the cumbersome iterations of the Newton-Raphson method.

\subsection{Optimization method}\label{optimize}
The network parameters can be optimized by minimizing the loss function (\ref{loss}) through the gradient descent method,
\begin{equation}\label{sgd}
    \theta_{k+1} = \theta_{k} - \alpha_k \nabla_\theta \mathcal{L}(\theta),\quad k=0,1,\cdots,K,
\end{equation}
where $\theta$ denotes parameters ${\theta_1, \theta_2}$ and $\theta_3$ (for nonlinear problem), $\alpha_k$ is the learning rate and $K$ is the total number of optimization steps. This requires calculating the gradient $\nabla_\theta L(\theta)$ which depends upon the solution of PDE (\ref{Learning-Correction}). The most straightforward method for computing the gradient is via the automatic differentiation of Tensorflow/Pytorch. Nevertheless, this technique incurs high computational and storage costs since the intermediate quantities of the iterations in our forward PDE solver, FEM, must be preserved to calculate the gradient. In the following, we derive the adjoint method for optimizing the linear multi-continuum model (\ref{Correction-linear}). \txrev{In the nonlinear problem with the forward solver PINN, the network parameters are optimized through \textbf{direct BP}. \textbf{Direct BP} denotes that the gradients of the loss function with respect to the network parameters are calculated by the automatic differentiation, as opposed to the adjoint method.}

\textbf{Optimization for linear equation.} The adjoint method \cite{bradley2013pde} provides a computationally efficient way to calculate gradients for PDE-constrained optimization problems. It requires solving the adjoint PDE per parameter update. The dimension of the adjoint PDE matches the dimension PDE (\ref{Learning-Correction}). The adjoint PDE of equation (\ref{Learning-Correction}) is given by
\begin{equation}\label{adjoint-equation}
\begin{aligned}
 & \frac{\partial \lambda_1}{\partial t} + \text{div}(\kappa_1^\star(\mathbf{x})\nabla \lambda_1) - \sigma(\mathbf{x}; \theta_2)(\lambda_1 - \lambda_2) = u_1(\mathbf{x},t) - U(\mathbf{x},t), \\
 & \frac{\partial \lambda_2}{\partial t} + \text{div}(\kappa_2(\mathbf{x}; \theta_1)\nabla \lambda_2) -\sigma(\mathbf{x}; \theta_2)(\lambda_2 - \lambda_1)= 0,
 \end{aligned}
\end{equation}
where $\lambda_1$ and $\lambda_2$ are the adjoint variable which are functions of coordinate $\mathbf{x}$ and time $t$. The initial and boundary conditions are given
$$ \lambda_1(\mathbf{x}, t=T) = \lambda_2(\mathbf{x}, t=T) = 0, \quad \text{in} \  \Omega \qquad
    \lambda_1(\mathbf{x}, t) = \lambda_2(\mathbf{x}, t) = 0, \quad \text{on} \  \partial\Omega\times [0,T]. $$
Then the gradient of the loss function (\ref{loss}) is
\begin{equation}\label{gradient}
\begin{aligned}
 \nabla_{\theta_1,\theta_2} \mathcal{L} =&  \int_0^T \int_{\Omega} -\lambda_2 \text{div}( \frac{\partial \kappa_2(\mathbf{x}; \theta_1)}{\partial \theta_1}\nabla u_2(\mathbf{x},t))d\mathbf{x}dt\\
 &+ \int_0^T \int_{\Omega} \lambda_1 \frac{\partial \sigma(\mathbf{x};\theta_2)}{\partial \theta_2} (u_1(\mathbf{x},t)-u_2(\mathbf{x},t)) + \lambda_2 \frac{\partial \sigma(\mathbf{x};\theta_2)}{\partial \theta_2} (u_2(\mathbf{x},t)-u_1(\mathbf{x},t))d\mathbf{x}dt.
\end{aligned}
\end{equation}
The derivation procedure of the adjoint equation (\ref{adjoint-equation}) is shown in \textbf{Appendix} A. \txrev{By this means, we calculate the derivative $\nabla_{\theta_1, \theta_2}\mathcal{L}$ by solving the adjoint PDE and calculating the derivatives $\frac{\partial \kappa_2(\mathbf{x};\theta_1)}{\partial \theta_1}$ and $\frac{\partial \sigma(\mathbf{x};\theta_2)}{\partial \theta_2}$ through the automatic differentiation.}

We calculate the numerical solution of adjoint equation (\ref{adjoint-equation}) also by FEM. The discretization is
\begin{equation}\label{FE-adjoint}
\begin{aligned}
 & M\frac{\partial \lambda_1}{\partial t} - A_1(\kappa_1^\star)\lambda_1 -\sigma(\mathbf{x};\theta_2)M(\lambda_1-\lambda_2) = M(u-U),\\
 & M\frac{\partial \lambda_2}{\partial t} - A_2(\kappa_2(\mathbf{x};\theta_1))\lambda_2 -\sigma(\mathbf{x};\theta_2)M(\lambda_2-\lambda_1) = 0,
\end{aligned}
\end{equation}
where $A_1,A_2$ are the stiffness matrices matrices, and $M$ is the mass matrix, as shown in Section \ref{forward}. The gradient (\ref{gradient}) can be reformulated as
\txrev{
\begin{equation}\label{FE-gradient}
\begin{aligned}
 \nabla_{\theta_1,\theta_2} \mathcal{L} =&  \int_T^0  \lambda_2^T \nabla_{\theta_1} A_2(\kappa_2)u_2dt\\
 +& \int_T^0 \lambda_1^T \nabla_{\theta_2} \sigma(\theta_2)M(u_1-u_2)dt\\
 +& \int_T^0 \lambda_2^T \nabla_{\theta_2} \sigma(\theta_2)M(u_2-u_1)dt.
\end{aligned}
\end{equation}}

\textbf{Optimization for nonlinear equation with PINN solver.} For the nonlinear multi-continuum model (\ref{Learning-Correction}), PINN is used as our PDE solver as discussed in section \ref{nonlinear}, it avoids the iterations when using FEM and the Newton-Raphson method. The network parameters $\theta_1, \theta_2$ and $\theta_3$ are optimized by direct BP in which the derivatives are calculated by automatic differentiation.

\section{Experiments}\label{experiments}
In this section, we will present some numerical experiments and demonstrate the performance of the proposed learning-based multi-continuum model. \txrev{Consider the flow equation on a unit square domain $\Omega = [0,1]\times[0,1]\subset \mathbb{R}^2$.} In the following experiments, we apply zero Dirichlet boundary conditions and zero initial conditions. The reference solutions are obtained by averaging the fine-scale solution over each coarse block. \txrev{The trusted data are sampled from the reference solutions. Further discussion on the trusted data can be found in \textbf{Appendix} B.}

We calculate the numerical solution of PDE (\ref{Learning-Correction}) by FEM. The loss function (\ref{loss}) is reformulated as
\begin{equation}\label{dis-loss}
\mathcal{L}=\frac{\sum\limits_{i,j=1}^{N_m} \sum\limits_{n=1}^{N_t} (u(\mathbf{x}_i, t_n) - U(\mathbf{x}_i, t_n))^TM_{ij}(u(\mathbf{x}_j, t_n) - U(\mathbf{x}_j, t_n))}{\sum\limits_{i,j=1}^{N_m} \sum\limits_{n=1}^{N_t} U(\mathbf{x}_i, t_n)^TM_{ij}U(\mathbf{x}_j, t_n)}\times 100\%,
\end{equation}
where $\{u(\mathbf{x}_i, t_n), U(\mathbf{x}_i,t_n)\}_{i=1,n=1}^{{N_m},N_t}$ are numerical solutions of our proposed model and the trusted data respectively, and $M_{ij}$ is corresponding element of the mass matrix.

\textbf{Model architecture.} The $\kappa$-net is employed to learn the permeability field of the second continuum. In the subsequent experiments, the $\kappa$-net is a multi-layer perceptron (MLP) with 5 hidden layers. \txrev{The number of hidden neurons is set to $100$ and the activation function is Tanh in each hidden layer.} In order to ensure that the permeability matrix $\kappa_2$ is a positive definite matrix for all coordinates and times, the $\kappa$-net $\kappa_2(\mathbf{x},t;\theta_1)$ is constrained to learn only the diagonal elements, with the additional requirement that they be greater than 0. Specifically, the learned permeability field $\kappa_2(\mathbf{x},t;\theta_1)=\left(\begin{array}{lc} \kappa_2^{11}(\mathbf{x},t;\theta_1) & 0 \\ 0 & \kappa_2^{22}(\mathbf{x},t;\theta_1) \end{array}\right)$, and the absolute value function is taken as the activation function in the output layer to constrain $\kappa_2^{11}(\mathbf{x},t;\theta_1) >0 ,\kappa_2^{22}(\mathbf{x},t;\theta_1) >0 $. 

The $\sigma$-net learns transfer coefficients for our learning-based multi-continuum model. \txrev{The $\sigma$-net used in our experiments is an MLP comprising $5$ hidden layers, with $LeakyReLU$ as the activation function and $100$ hidden neurons per hidden layer.} The $LeakyReLU$ function is defined as
$$ LeakyReLU(x)=\left\{\begin{aligned} &x, & x>0,\\ &cx, & x<=0. \end{aligned} \right.$$
We set $c=0.2$ in all experiments.

In our experiment of nonlinear flow equations, the PINN solver is adopted and parameterized as an MLP with 7 hidden layers. The number of hidden neurons is 100, and the activation function is Tanh. Furthermore, the homogenized permeability $\kappa_1^\star$ represents the average value of each sub-domain in $\Omega$. To calculate the homogenized permeability for each grid point, a convolution operation is employed to act as an interpolation function, with a convolution kernel size of $2\times2$. \txrev{The elements of this kernel are fixed as a constant value of 1/4, without training}.

\textbf{Hyperparameters.} The spatial partition $\Gamma_h$ is the uniform mesh of the computational domain $\Omega$, and the mesh size is set $h=0.1$. The number of linear basis functions used in FEM is the same as the number of mesh points. We optimize the network parameters by Adam algorithm \cite{kingma2014adam}, and the learning rate is set to $1\times10^{-4}$. The total simulation time $T$ and the time step size $\tau$ are separately presented in each of subsections 4.1-4.3.

\begin{figure}[h]
        \centering
        \subfigure[Permeability field] {\includegraphics[width=0.23\columnwidth]{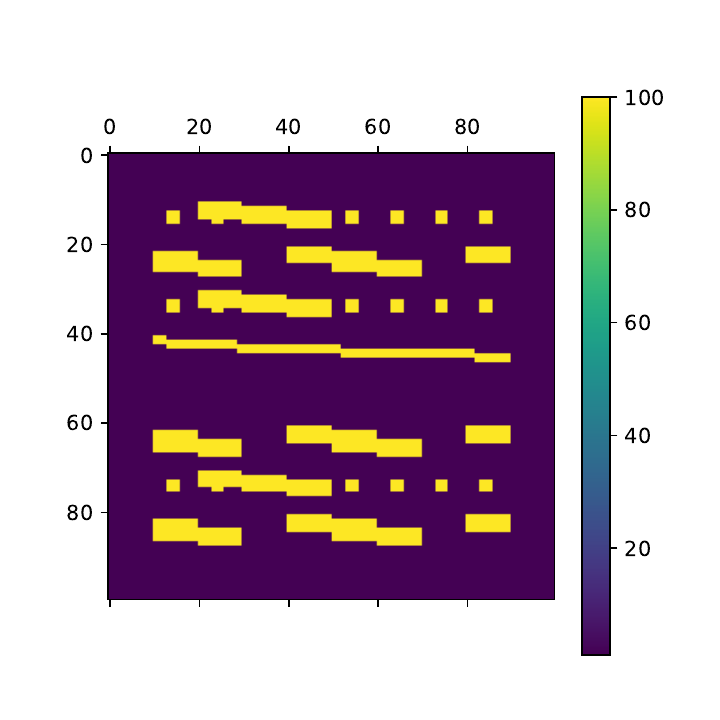}}
        \subfigure[Source term] {\includegraphics[width=0.24\columnwidth]{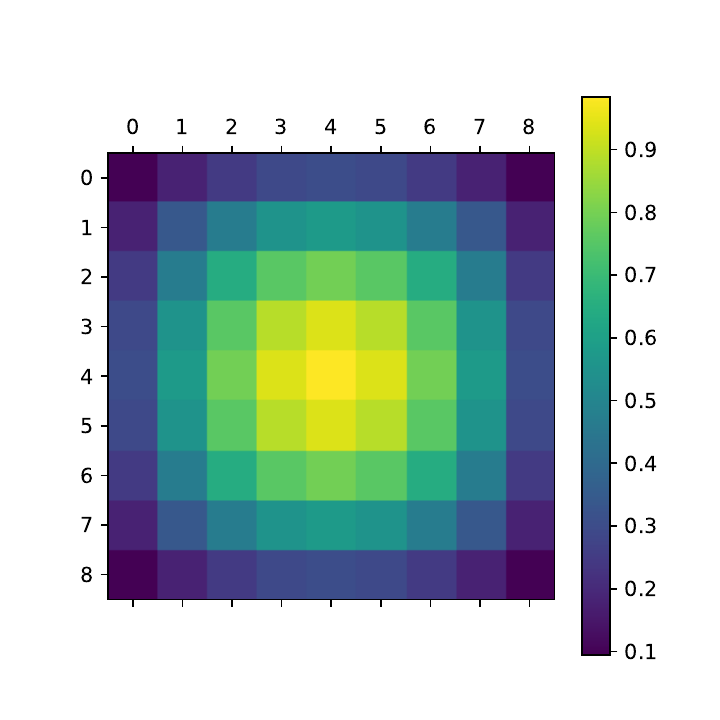}}
        \subfigure[Homogenized $\kappa^\star_{11}$] {\includegraphics[width=0.24\columnwidth]{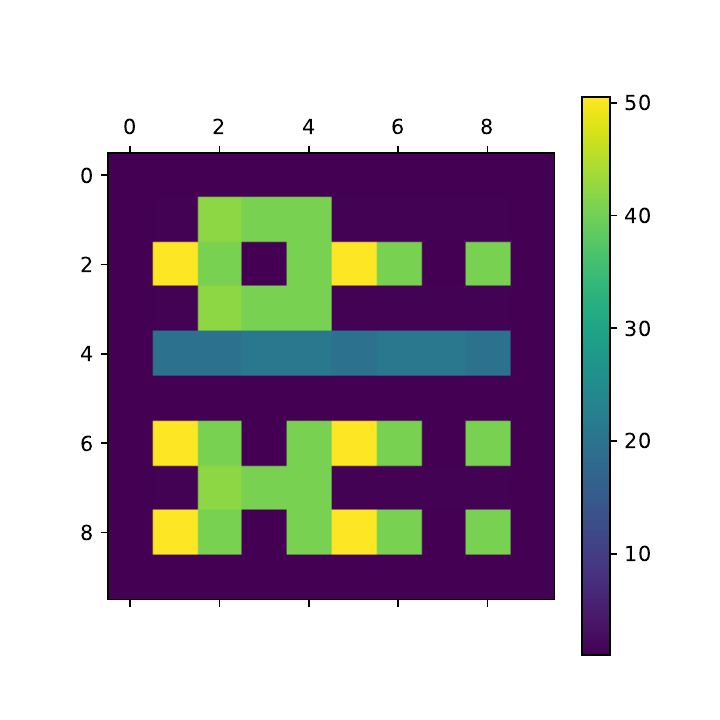}}
        \subfigure[Homogenized $\kappa^\star_{22}$] {\includegraphics[width=0.24\columnwidth]{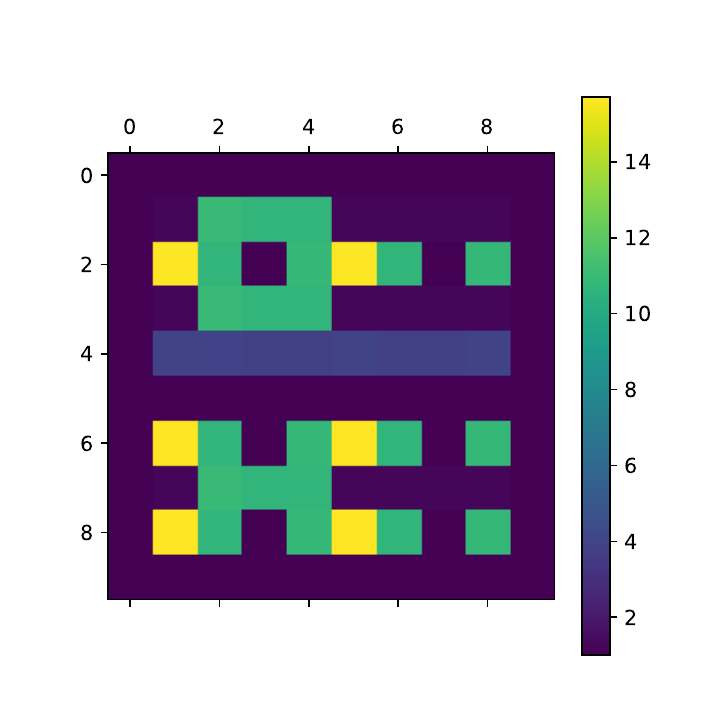}}
        \caption{The configuration of the permeability field and the source term of the linear flow equation. (a) The permeability field on the fine grid. (b) The source term. (c) The value $\kappa^\star_{11}$ of the homogenized permeability. (d) The value $\kappa^\star_{22}$ of the homogenized permeability.}
        \label{linear-permeability}
\end{figure}

\begin{figure}[t]
        \centering
        \subfigure[Direct back-propagation]{\includegraphics[width=0.3\columnwidth]{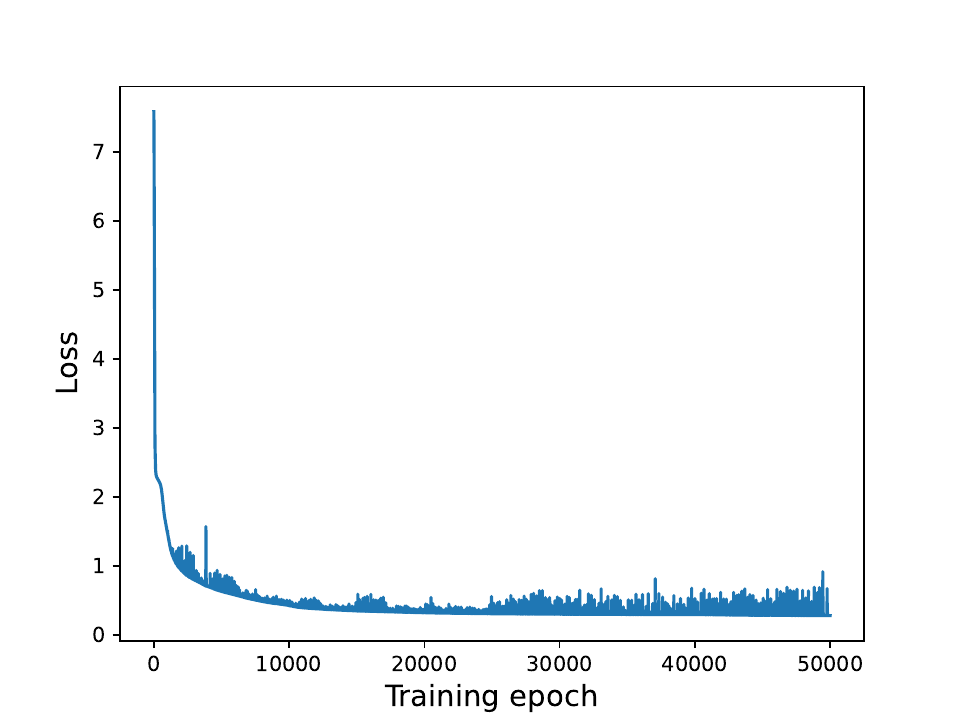}}
        \hspace{.25in}
        \subfigure[Adjoint method]{\includegraphics[width=0.3\columnwidth]{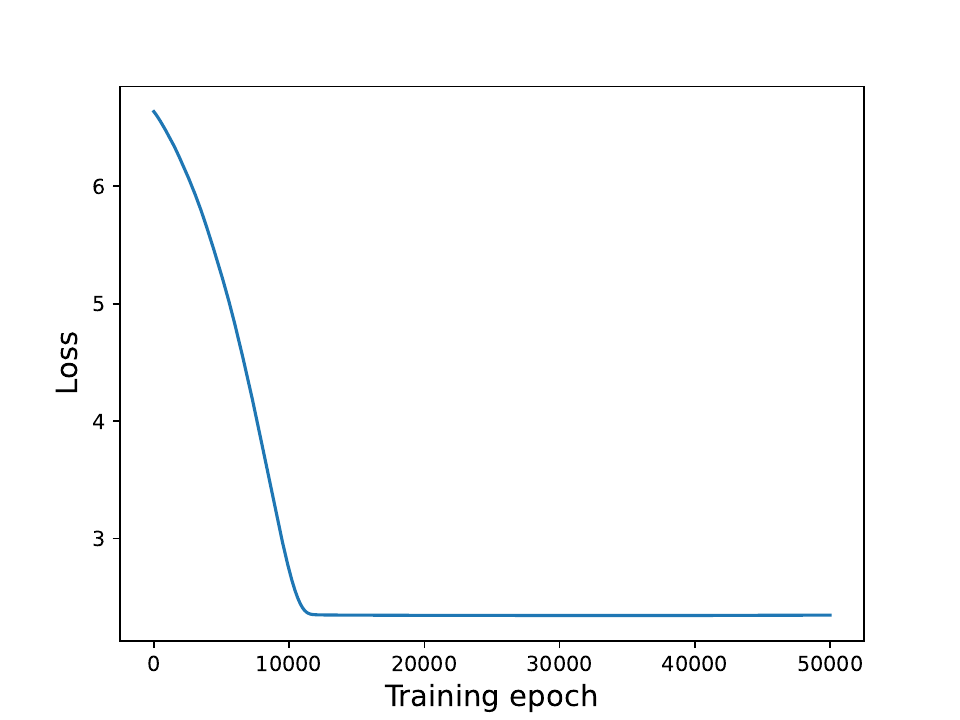}}
        \caption{The training trajectories of the loss function. The gradients of netwrok parameters are calculated by (a) direct BP and (b) the adjoint method. }
        \label{linear-loss}
\end{figure}

\subsection{Example 1: Linear flow equation}
In the first example, the configuration of the permeability field and the source term $f(\mathbf{x})$ are shown in \textbf{Figure} \ref{linear-permeability} (a-b). The value of permeability is $100$ in the channel and $1$ in the background. The homogenized equation is obtained by numerical homogenization, where the homogenized permeability $\kappa_1^\star=\left(\begin{array}{lc} \kappa^\star_{11} & \kappa^\star_{12} \\ \kappa^\star_{21} & \kappa^\star_{22} \end{array}\right)$. The diagonal elements $\kappa^\star_{11}$ and $\kappa^\star_{22}$ are shown in \textbf{Figure} \ref{linear-permeability} (c-d), which are globally consistent with the configuration of the original permeability field. The total simulation time $T=1$ and the time step size $\tau=0.1$.

The gradients of the loss function (\ref{dis-loss}) with respect to the network parameters are calculated using both direct BP and the adjoint method to optimize the proposed learning-based multi-continuum model (\ref{Correction-linear}). \textbf{Figure} \ref{linear-loss} shows the training trajectories of the loss function under two gradient calculation methods. It is evident that the loss function decreases rapidly and saturates after 12000 epochs for both methods. While direct BP can achieve a loss function of $O(10^{-1})$, the adjoint method can reach a loss function of $O(1)$. %This might be due to the difference in the type of gradient descent.

\begin{figure}[t]
        \centering
        \subfigure[Reference solution $t=0.1$]
        {\includegraphics[width=0.25\columnwidth]{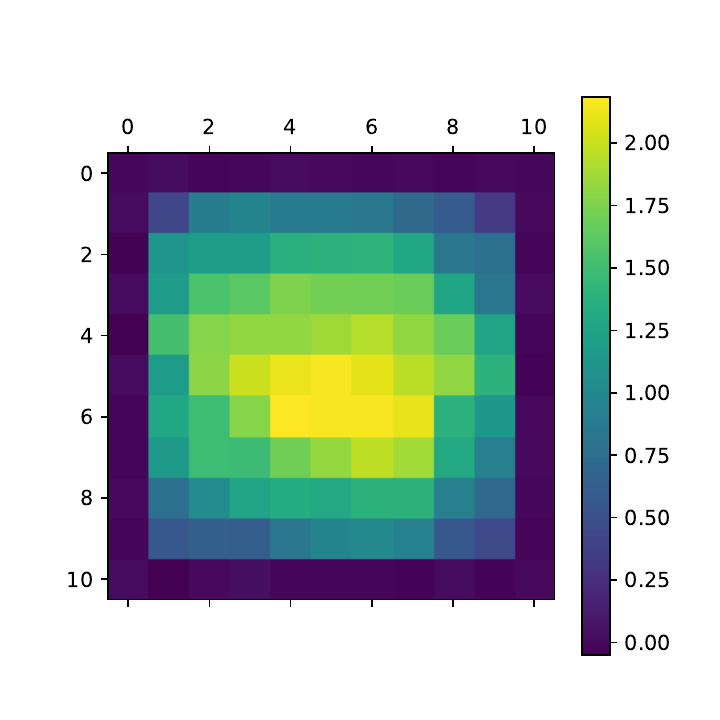}}
        \subfigure[Homogenized solution $t=0.1$] {\includegraphics[width=0.25\columnwidth]{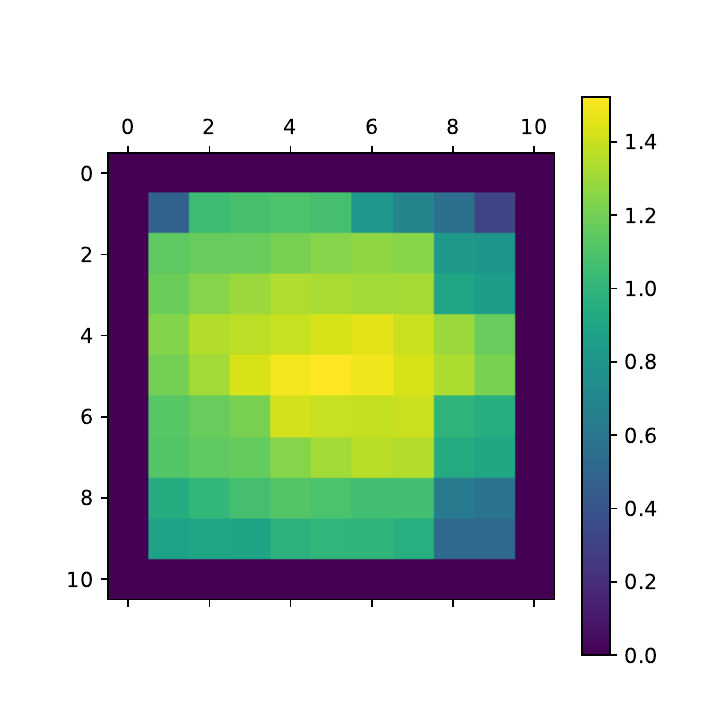}}
        \subfigure[Corrected solution $t=0.1$]
        {\includegraphics[width=0.25\columnwidth]{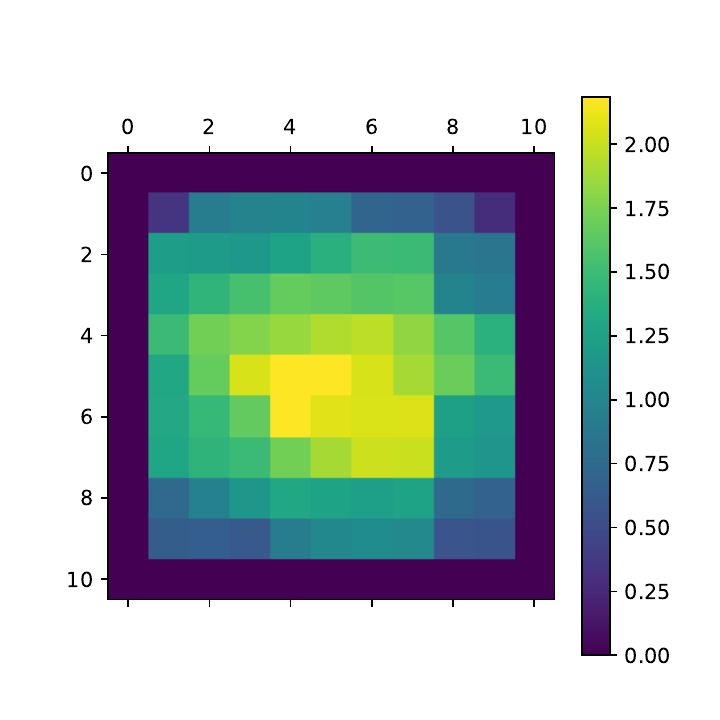}}\\

        \subfigure[Reference solution $t=1$]
        {\includegraphics[width=0.25\columnwidth]{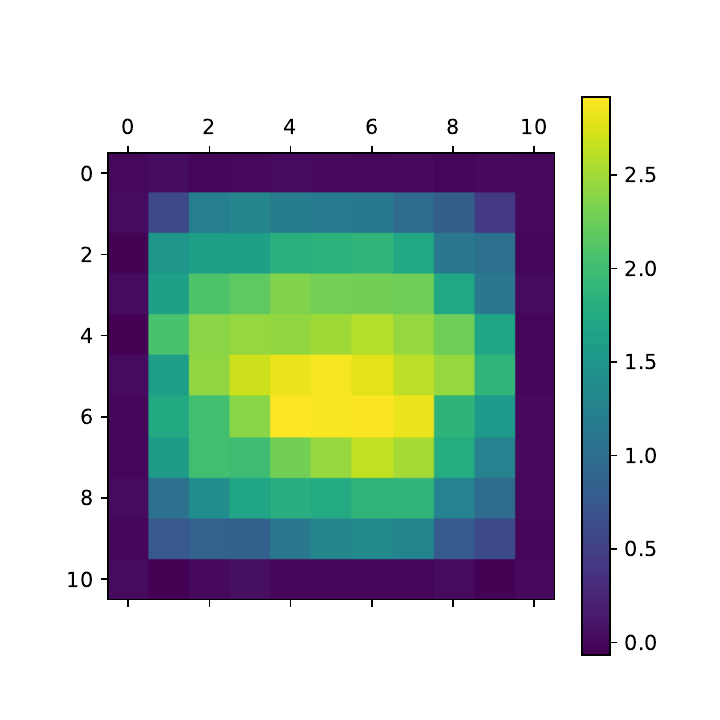}}
        \subfigure[Homogenized solution $t=1$] {\includegraphics[width=0.25\columnwidth]{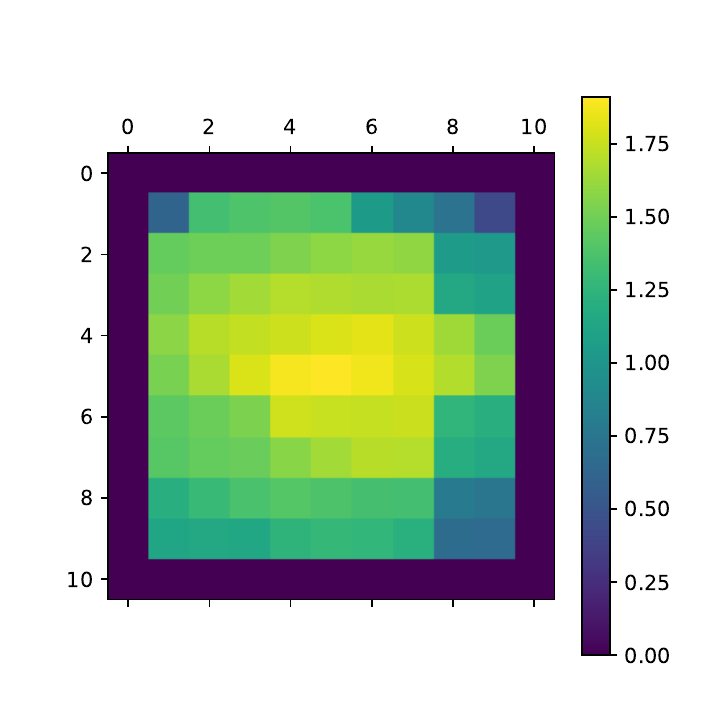}}
        \subfigure[Corrected solution $t=1$]
        {\includegraphics[width=0.25\columnwidth]{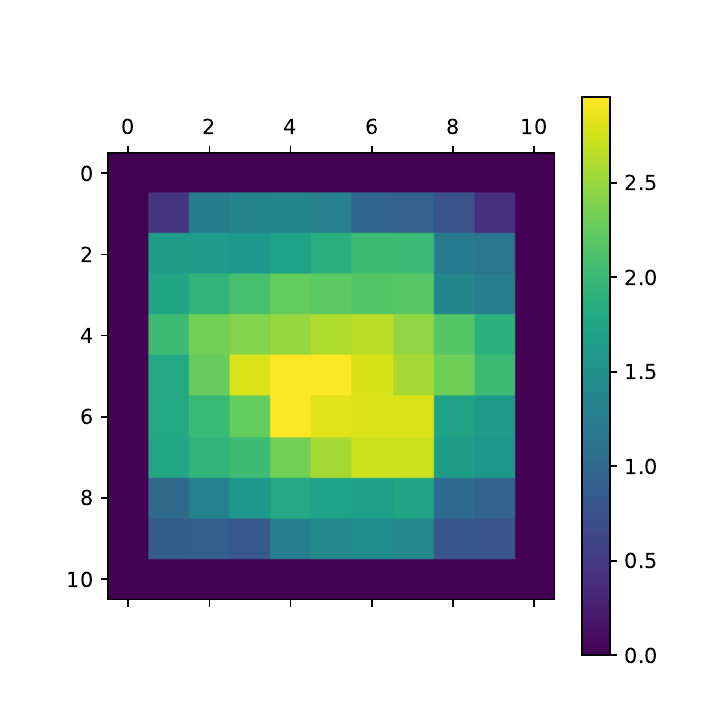}}

        \caption{The FEM solutions of the linear flow equation, the homogenized equation and our proposed learning-based multi-continuum equation. (a)(d) The reference solutions at $t=0.1$ and 1. (b)(e) The homogenized solutions at $t=0.1$ and 1. (c)(f) The solutions of our proposed multi-continuum equation at $t=0.1$ and 1.}
        \label{linear-solution}
\end{figure}

\textbf{Table}~\ref{tbl1} lists the relative $L_2$ error of the homogenized solutions and corrected solutions of our learning-based multi-continuum model with respect to the reference solutions. It can be seen that the numerical homogenization results in a significant relative $L_2$ error for each time step, which accumulates over time. The relative $L_2$ error increases from 5.7956\% at time $t=0.1$ to 7.7845\% at time $t=1$. However, our proposed learning-based multi-continuum model corrects the homogenized equation resulting in more accurate numerical solutions. The relative $L_2$ error of our multi-continuum model optimized by direct BP decreased from $0.1162\%$ ($t=0.1$) to $0.1084\%$ ($t=1$). The multi-continuum model optimized by the adjoint method for gradient calculation performs relatively poorly, with a relative $L_2$ error of approximately 2.33\% for all time steps. \textbf{Figure} \ref{linear-solution} further displays the reference solutions, homogenized solutions, and the multi-continuum model's solutions at different time steps. Note that all of these solutions are FEM solutions. \txrev{The figure reflects that our proposed learning-based multi-continuum model corrects the homogenized equation and generates the numerical solutions that closely match the reference solutions.}

\begin{figure}[h]
        \centering
        \subfigure[Learned $\kappa_2^{11}$]
        {\includegraphics[width=0.25\columnwidth]{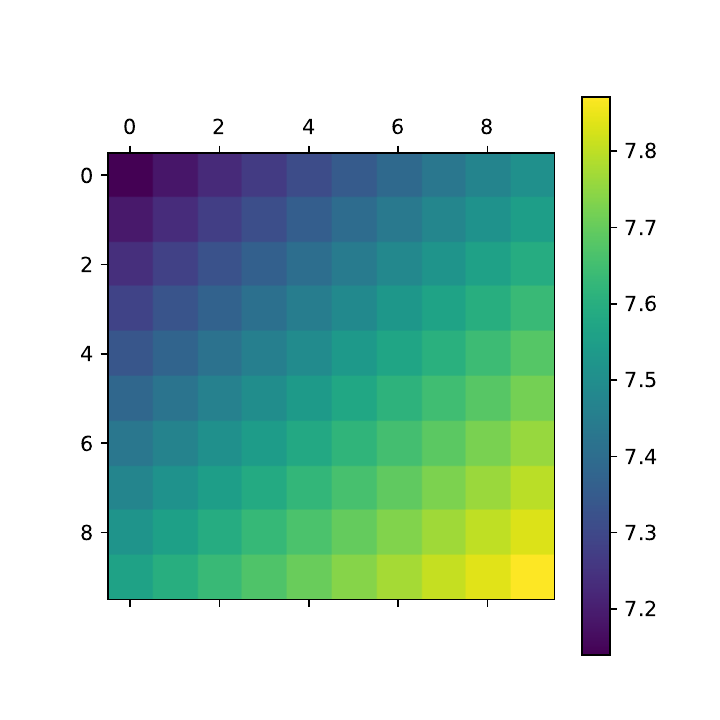}}
        \subfigure[Learned $\kappa_2^{22}$]
        {\includegraphics[width=0.25\columnwidth]{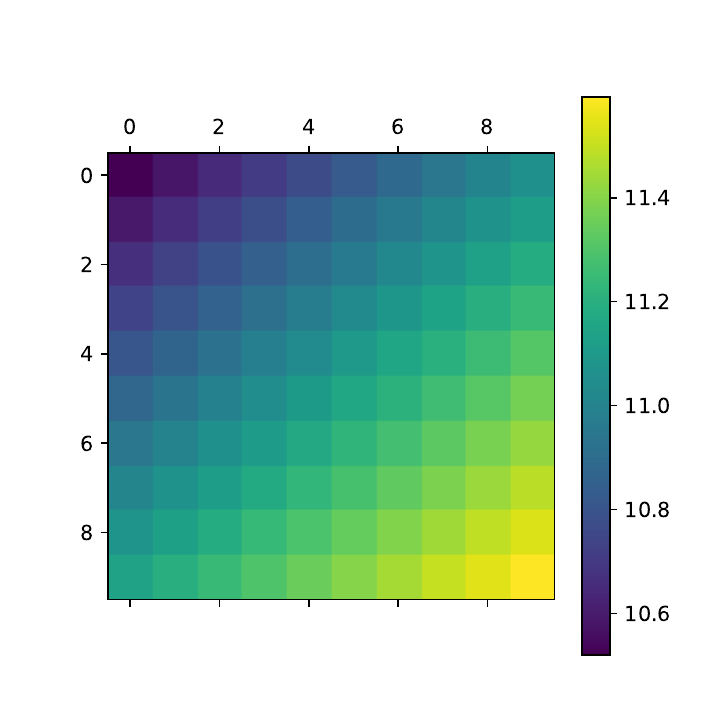}}
        \subfigure[Learned $\sigma$]
        {\includegraphics[width=0.25\columnwidth]{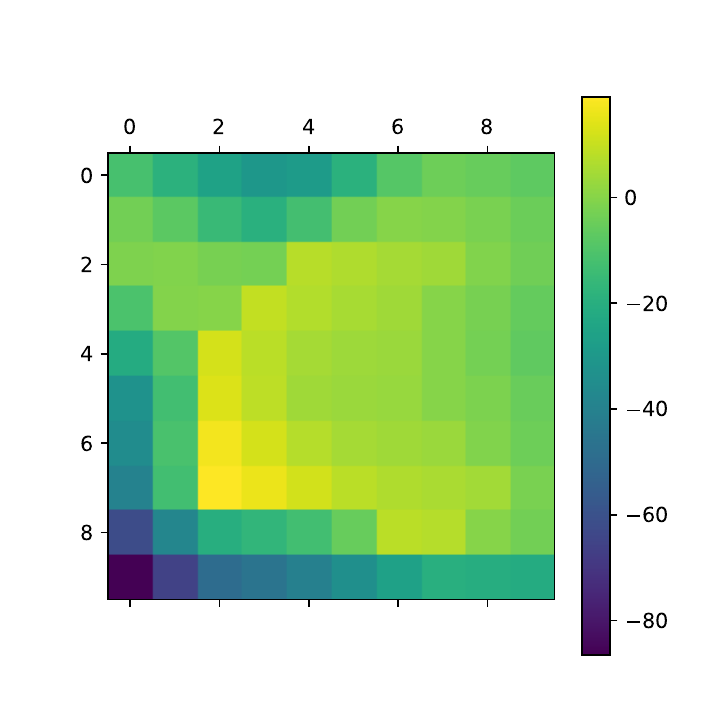}}

        \caption{The learned permeability field $\kappa_2$ and transfer coefficient $\sigma$. (a-b) The diagonal elements $\kappa_2^{11}$ and $\kappa_2^{22}$ in domain $\Omega$. (c) The learned transfer coefficient $\sigma$ in domain $\Omega$.}
        \label{linear-learned-parameter}
\end{figure}

\textbf{Figure}~\ref{linear-learned-parameter} visualises the values of learned permeability $\kappa_2(\mathbf{x};\theta_1)$ and the learned transfer coefficient $\sigma(\mathbf{x};\theta_2)$. As discussed in \textbf{model architecture} of Section \ref{experiments}, the learned permeability $\kappa_2(\mathbf{x}, \theta_1)=\left(\begin{array}{lc} \kappa_2^{11}(\mathbf{x}, \theta_1) & 0 \\ 0 & \kappa_2^{22}(\mathbf{x}, \theta_1) \end{array}\right)$. \textbf{Figure}~\ref{linear-learned-parameter} (a-b) shows the diagonal elements  $\kappa_2^{11}(\mathbf{x}, \theta_1)$ and $\kappa_2^{22}(\mathbf{x}, \theta_1)$. It can be seen that these two values have the same trend in the domain $\Omega$, gradually increasing from the upper left corner to the lower right corner, but with relatively small increases. As illustrated in Figure Figure \ref{linear-permeability} (c-d), the homogenized permeability field retains some average information of the original permeability field (Figure \ref{linear-permeability} (a)), while the learned permeability field $\kappa_2$ (Figure \ref{linear-learned-parameter} (a-b)) represents some missing or overestimate information from $\kappa^\star_1$. These two flow equations are composited into a new multi-continuum model in which the background and the fracture are written separately with some interaction. This multi-continuum model more accurately characterizes the complex dynamics within fractured media. \textbf{Figure}~\ref{linear-learned-parameter} (c) shows the learned transfer coefficient $\sigma(\theta_2)$. The absolute values of $\sigma$ are relatively small at the center of $\Omega$, which reflects that the correction of the homogenized equation by the second continuum is relatively weaker in the center region of $\Omega$. In other words, the numerical homogenization for this linear problem can portray the variation of the flow in the central region.

\begin{table}[htp]
\caption{The relative $L_2$ error (in percentage) of the homogenized solutions and solutions of the proposed multi-continuum model with respect to the reference solutions.}\label{tbl1}
\begin{tabular*}{\tblwidth}{@{}LLLLLLLLLLL@{}}
\toprule
 \diagbox{Method}{Time} & 0.1 & 0.2 & 0.3 & 0.4 & 0.5 & 0.6 & 0.7 & 0.8 & 0.9 & 1.0 \\ % Table header row
\midrule
 Homogenization &  5.7956 & 7.0050 & 7.5187 & 7.7016 & 7.7560 & 7.7774 & 7.7824 & 7.7838 & 7.7841 & 7.7845\\
 \makecell[l]{Multi-continuum model \\ (direct BP)}&  0.1162 & 0.1148 & 0.1106 & 0.1091 & 0.1086 & 0.1085 & 0.1084 & 0.1084 & 0.1084 & 0.1084\\
 \makecell[l]{Multi-continuum model \\ (adjoint method)} &  2.3355 & 2.3451 & 2.3379 & 2.3338 & 2.3321 & 2.3314 & 2.3311 & 2.3310 & 2.3310 & 2.3310\\
\bottomrule
\end{tabular*}
\end{table}

\subsection{Example 2: Linear flow equation with complex permeability field}
\begin{figure}[t]
        \centering
        \subfigure[Permeability field] {\includegraphics[width=0.24\columnwidth]{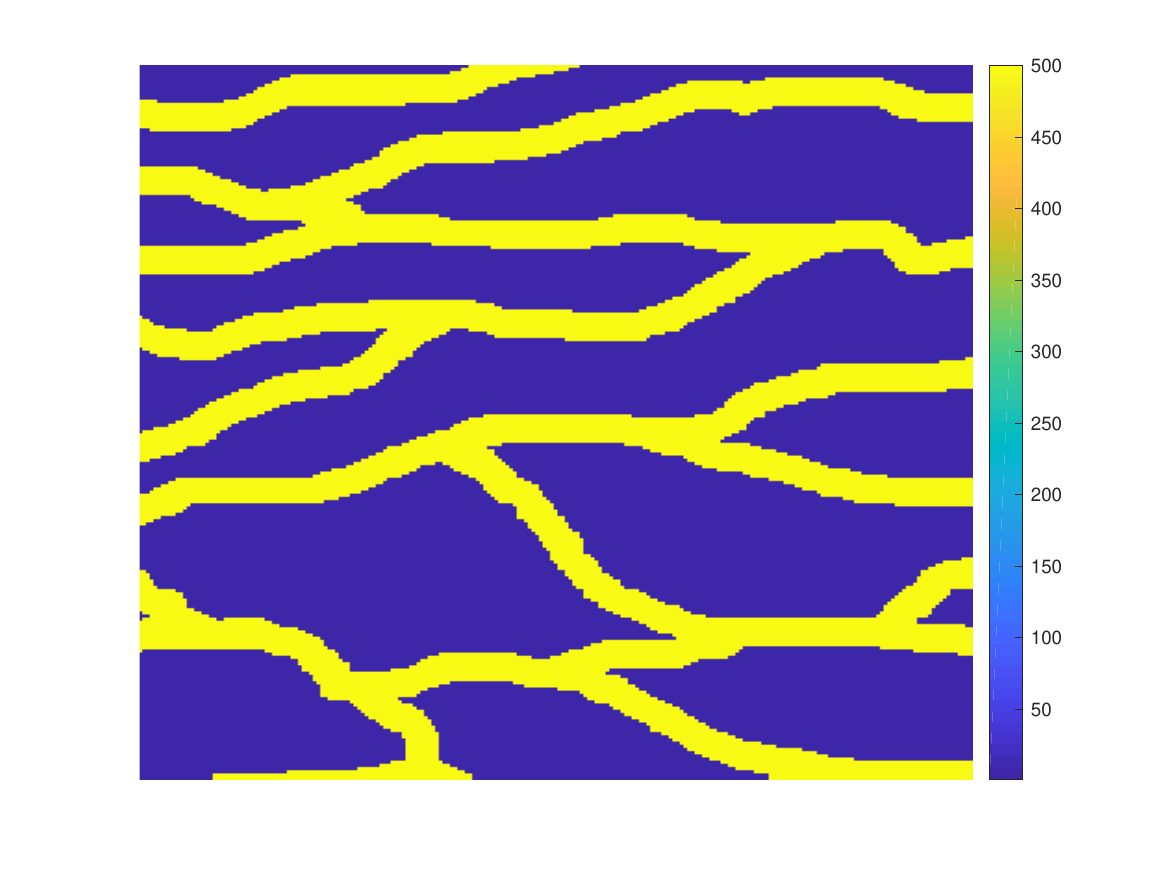}}
        \hspace{.35in}
        \subfigure[Source Term] {\includegraphics[width=0.24\columnwidth]{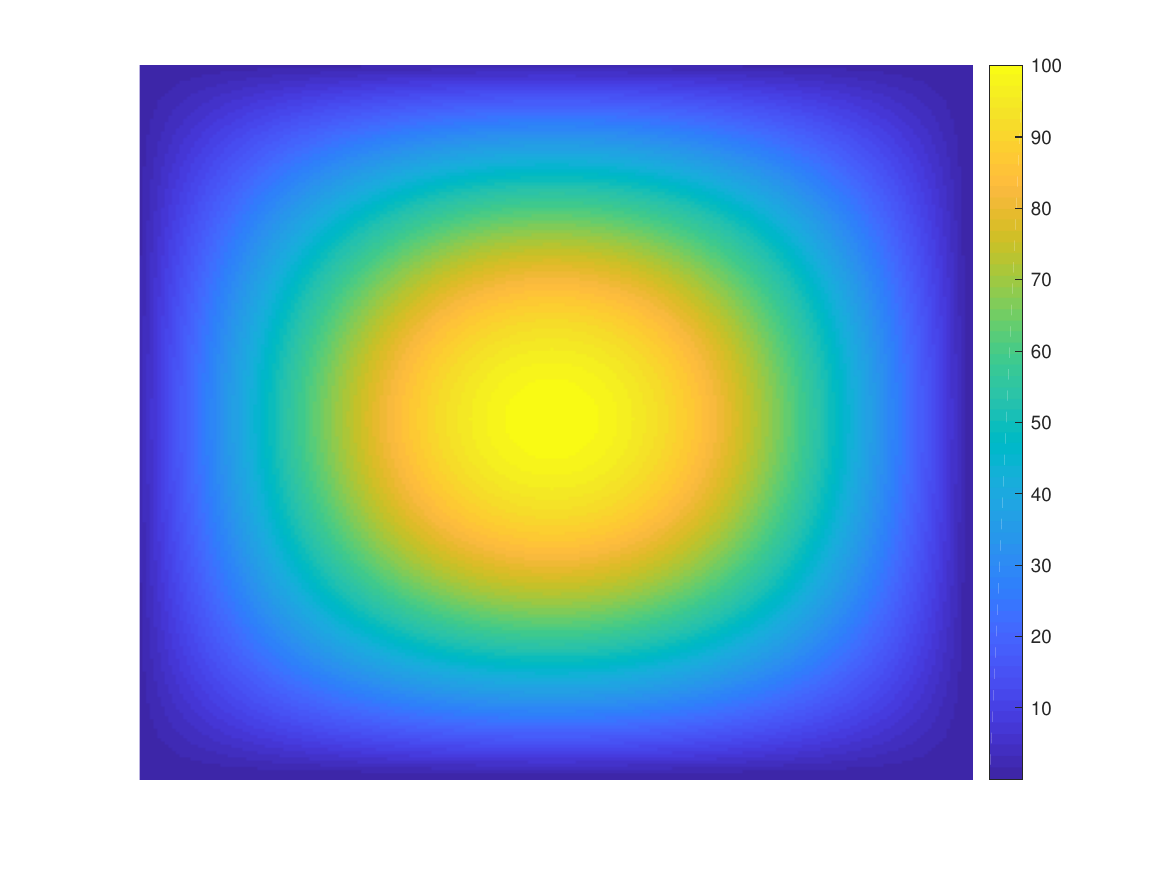}}
        \caption{The configuration of the permeability field and the source term. (a) The permeability field. (b) The source term.}
        \label{linear2-permeability}
    \end{figure}

The configuration of the permeability field and the source term $f(\mathbf{x})$ for this example are shown in \textbf{Figure} \ref{linear2-permeability}. This permeability field is closer to the real underground scene and is more complex. The value of permeability is $500$ in the channel and $1$ in the background. The homogenized equation is obtained by numerical homogenization. The total simulation time $T=0.001$ and the time step size $\tau=0.0001$.

\begin{table}[h]
\caption{The relative $L_2$ error (in percentage) of the homogenized solutions and solutions of the proposed multi-continuum model with respect to the reference solutions.}\label{tbl2}
\begin{tabular*}{\tblwidth}{@{}LLLLLLLLLLL@{}}
\toprule
 \diagbox{Method}{Time} & 0.0001 & 0.0002 & 0.0003 & 0.0004 & 0.0005 & 0.0006 & 0.0007 & 0.0008 & 0.0009 & 0.001 \\ % Table header row
\midrule
 Homogenization &  1.2710 & 2.2432 & 3.3052 & 4.4189 & 5.5602 & 6.7125 & 7.8634 & 9.0037 & 10.1263 & 11.2255\\
 Multi-continuum model &  0.5731 & 0.7715 & 0.9768 & 1.1847 & 1.3851 & 1.5715 & 1.7423 & 1.8996 & 2.0491 & 2.1995\\
\bottomrule
\end{tabular*}
\end{table}

\begin{figure}[t]
        \centering
        \subfigure[Reference solution $t=0.0001$]
        {\includegraphics[width=0.25\columnwidth]{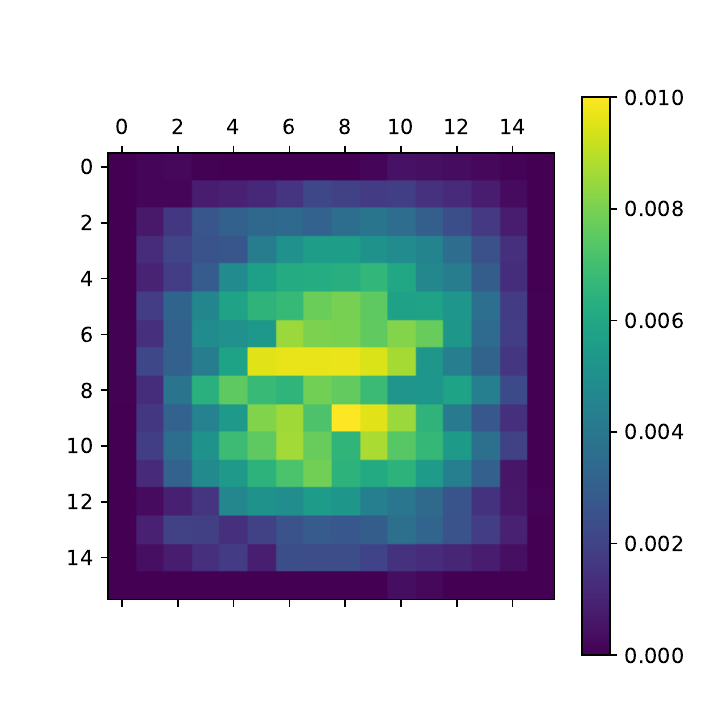}}
        \subfigure[Homogenized solution $t=0.0001$] {\includegraphics[width=0.25\columnwidth]{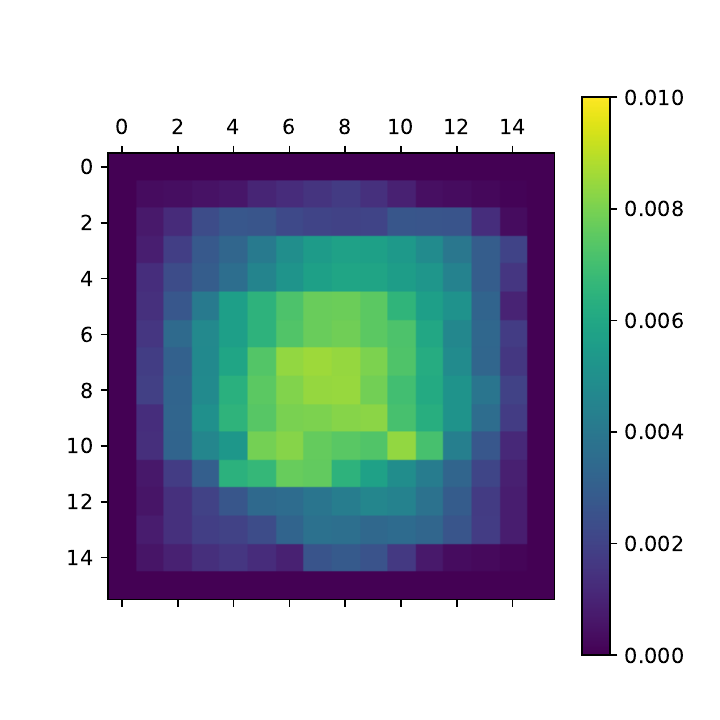}}
        \subfigure[Corrected solution $t=0.0001$]
        {\includegraphics[width=0.25\columnwidth]{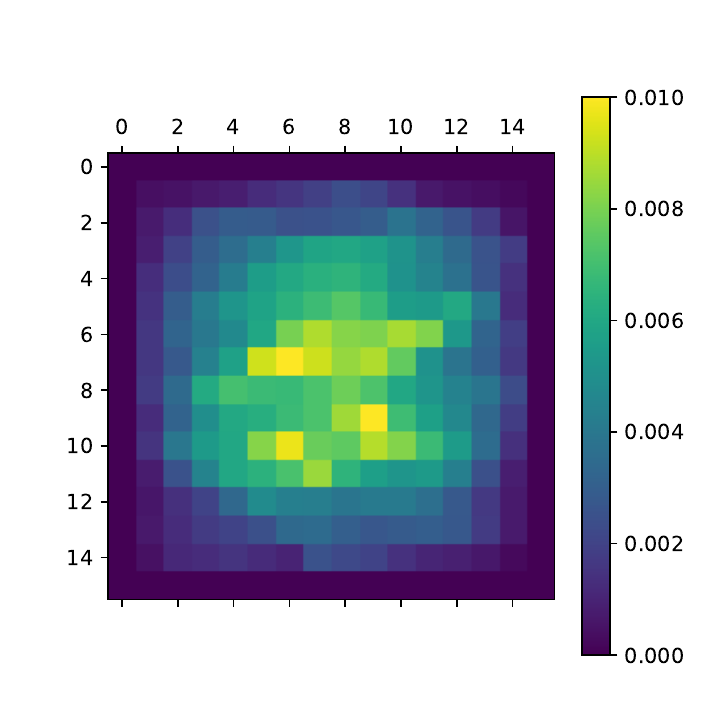}}\\

        \subfigure[Reference solution $t=0.0005$]
        {\includegraphics[width=0.25\columnwidth]{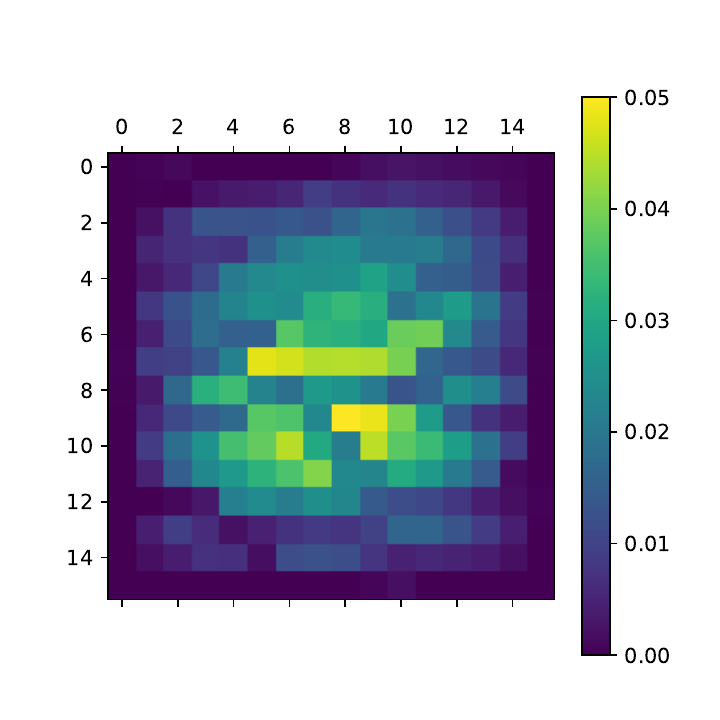}}
        \subfigure[Homogenized solution $t=0.0005$] {\includegraphics[width=0.25\columnwidth]{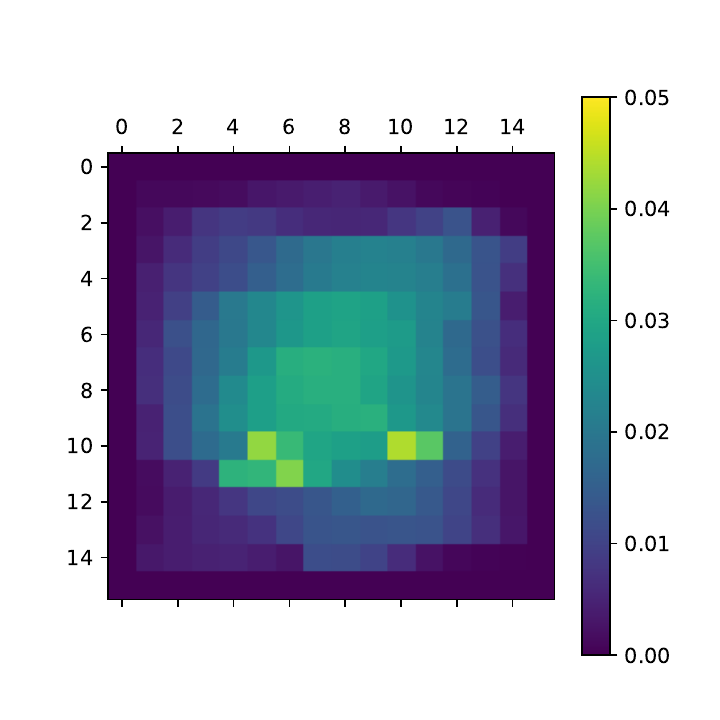}}
        \subfigure[Corrected solution $t=0.0005$]
        {\includegraphics[width=0.25\columnwidth]{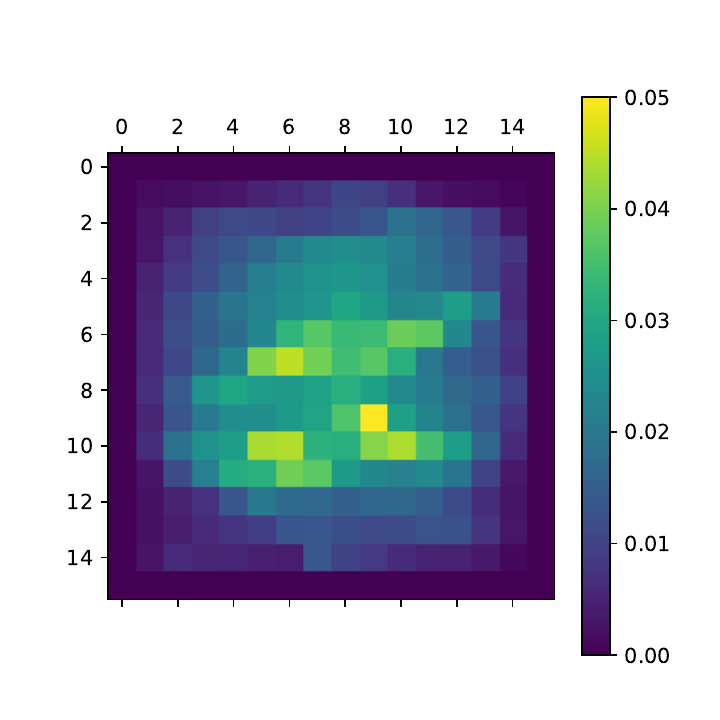}}

        \subfigure[Reference solution $t=0.001$]
        {\includegraphics[width=0.25\columnwidth]{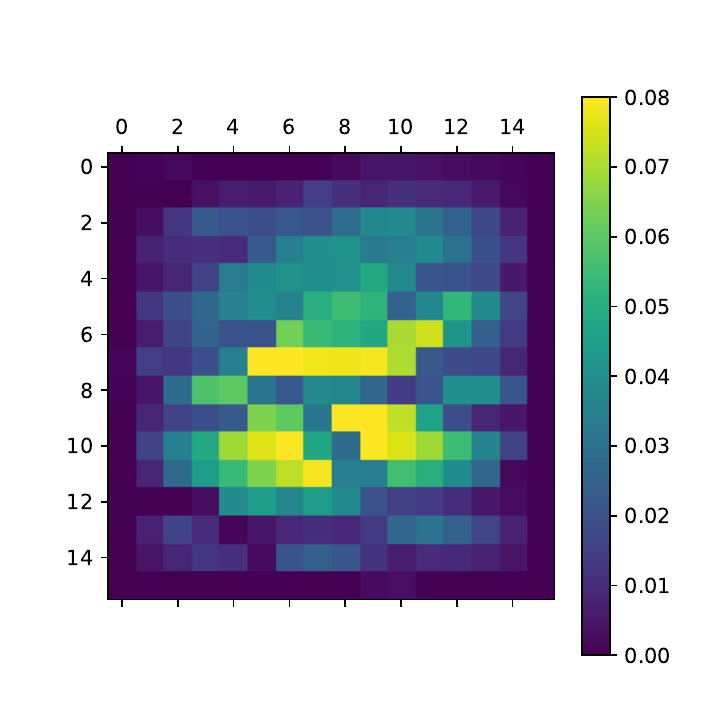}}
        \subfigure[Homogenized solution $t=0.001$] {\includegraphics[width=0.25\columnwidth]{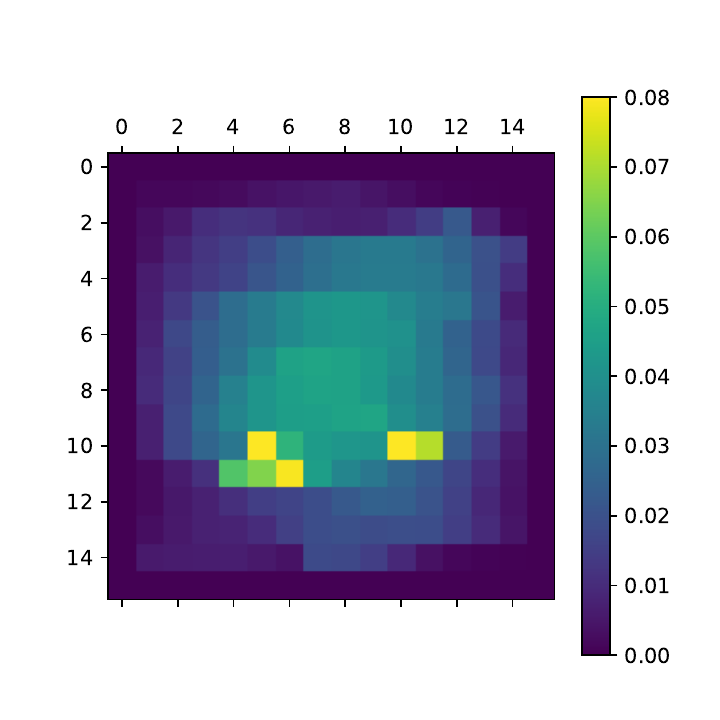}}
        \subfigure[Corrected solution $t=0.001$]
        {\includegraphics[width=0.25\columnwidth]{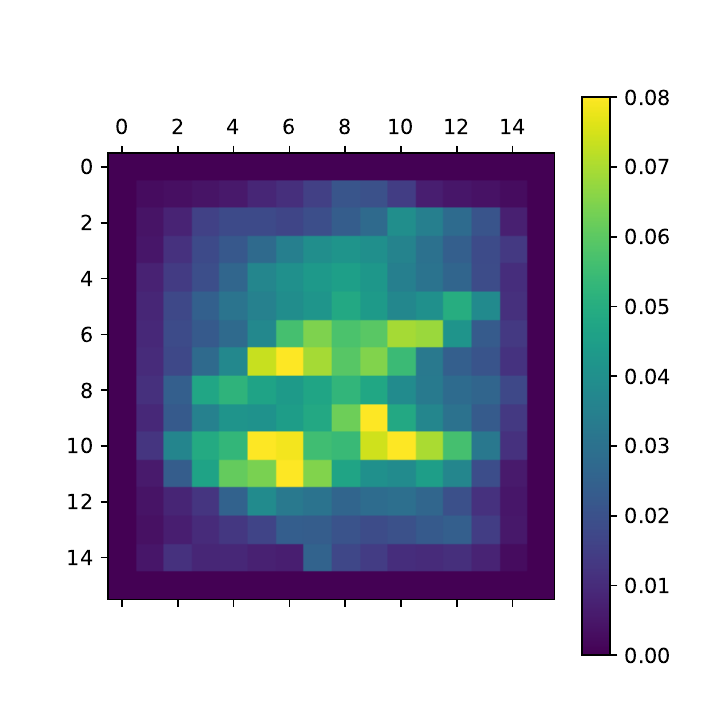}}

        \caption{The FEM solutions of the flow equation with complex permeability field, the homogenized equation and our proposed learning-based multi-continuum equation at time $t=0.0001,0.0005,0.001$. (a)(d)(g) The reference solutions. (b)(e)(h) The homogenized solutions. (c)(f)(i) The corrected solutions.}
        \label{linear2-solution}
\end{figure}

\begin{figure}[t]
        \centering
        \subfigure[Learned $\kappa_2^{11}$]
        {\includegraphics[width=0.25\columnwidth]{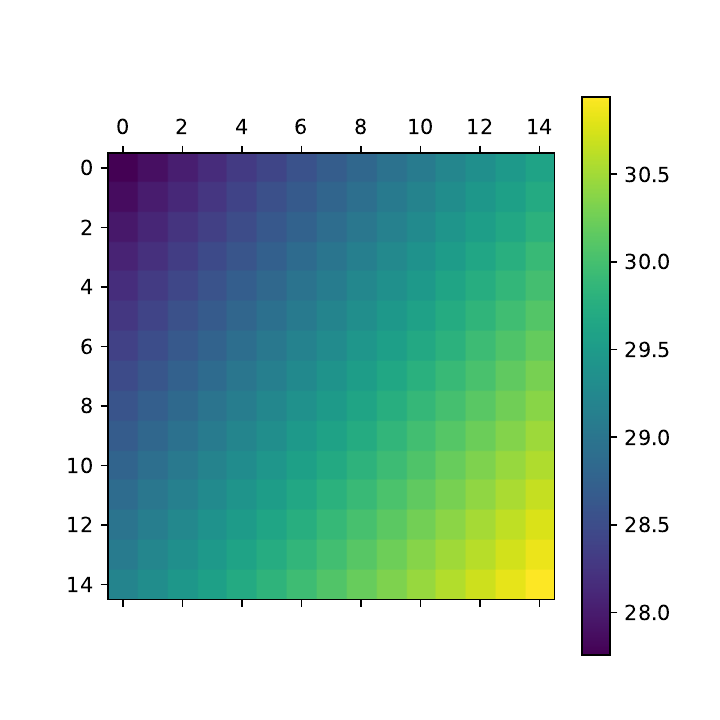}}
        \subfigure[Learned $\kappa_2^{22}$]
        {\includegraphics[width=0.25\columnwidth]{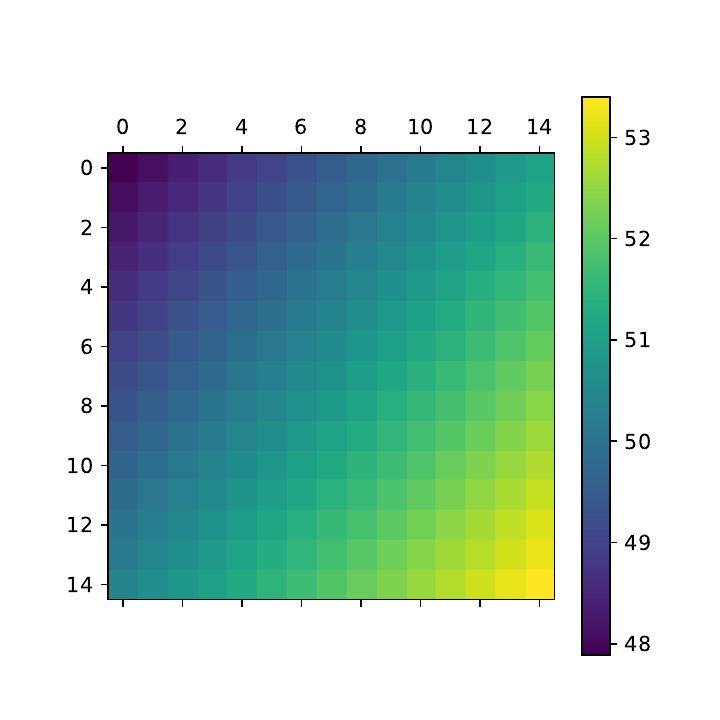}}
        \subfigure[Learned $\sigma$]
        {\includegraphics[width=0.25\columnwidth]{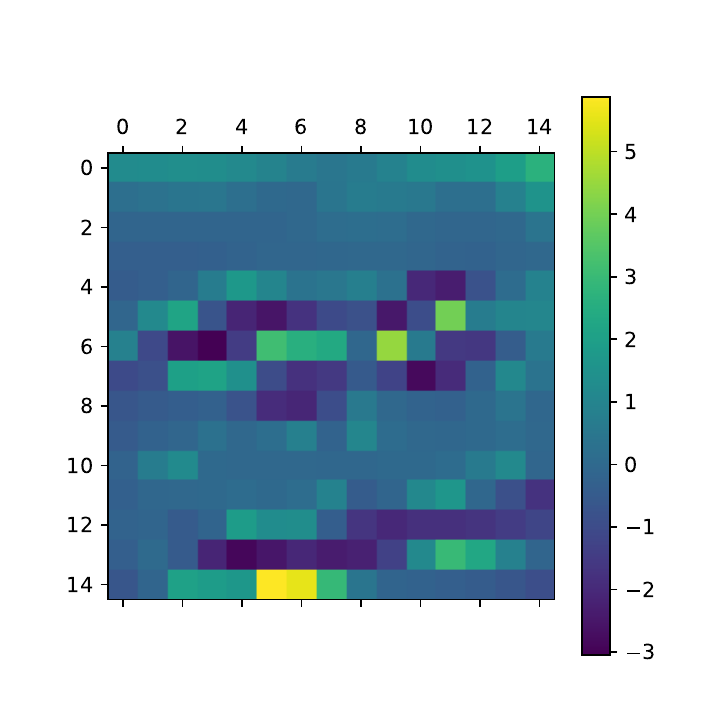}}
        \caption{The learned permeability field $\kappa_2$ and transfer coefficient $\sigma$. (a-b) The diagonal elements $\kappa_2^{11}$ and $\kappa_2^{22}$ in domain $\Omega$. (c) The learned $\sigma$.}
        \label{linear2-learned-parameter}
\end{figure}

\textbf{Table}~\ref{tbl2} lists the relative $L_2$ error of the homogenized solutions and the corrected solutions of our multi-continuum model with respect to the reference solutions. The relative $L_2$ error of the homogenized solutions of this problem gradually increases over time, from 1.2710\% at $t=0.0001$ to 11.2255\% at $t=0.001$. Our learning-based multi-continuum model can correct the homogenized equation resulting in a small relative $L_2$ error. At time $t=0.001$, the relative $L_2$ error of our multi-continuum model's solution is 2.1995\%, which is one-fifth of the relative $L_2$ error for numerical homogenization. Meanwhile, the results reflect that our learning-based multi-continuum model reduces the accumulation of errors in the numerical solutions over time. \textbf{Figure} \ref{linear2-solution} illustrates the numerical solutions of the original multiscale flow equation, the homogenized equation, and our proposed multi-continuum model at different time steps. The results demonstrate that our learning-based multi-continuum model effectively corrects the homogenized solutions and can capture the dynamics of the flow equation with this complex permeability field.

\textbf{Figure}~\ref{linear2-learned-parameter} illustrates the values of the learned permeability $\kappa_2(\mathbf{x};\theta_1)$ and the learned transfer coefficient $\sigma(\mathbf{x};\theta_2)$. The diagonal elements $\kappa_2^{11}(\mathbf{x};\theta_1)$ and $\kappa_2^{22}(\mathbf{x};\theta_1)$ of the learned permeability are shown in \textbf{Figure}~\ref{linear2-learned-parameter} (a-b). It is apparent that these two values vary throughout the whole domain $\Omega$, but the variation is relatively small. Similar to the first example, the learned permeability can be served to find the missing/over-estimated information for the numerical homogenization. \textbf{Figure}~\ref{linear2-learned-parameter} (c) displays the learned transfer coefficient $\sigma(\mathbf{x};\theta_2)$ throughout the entire domain $\Omega$. The result shows that the region of significant sigma fluctuations corresponds to the area of fluid flow alteration, while changes in the peripheral region are gradual, and the values are relatively small. These values of the learned transfer coefficient $\sigma$ reflect the impact of adding the second continuum to the numerical homogenization. 

\subsection{Example 3: Nonlinear flow equation with discontinuous source term}
\begin{figure}[t]
        \centering
        \subfigure[Permeability field] {\includegraphics[width=0.24\columnwidth]{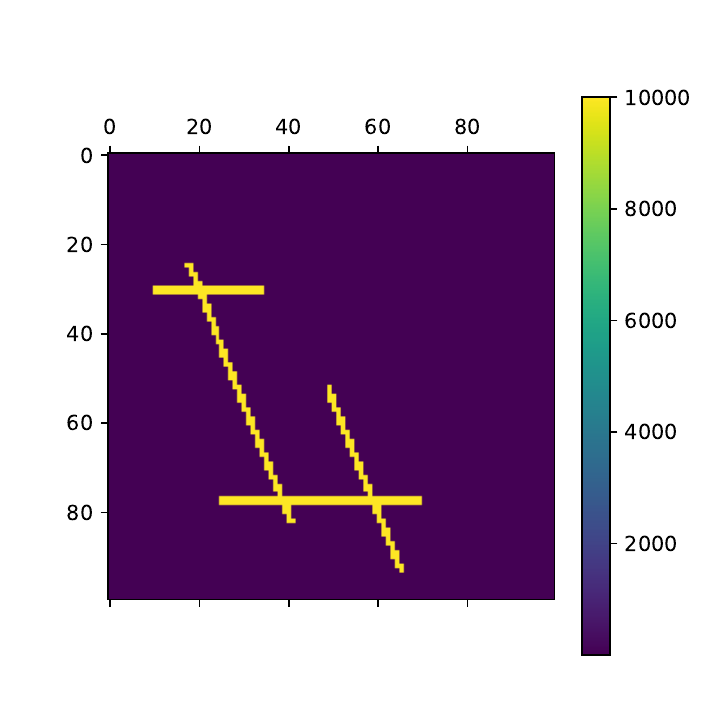}}
        \subfigure[Source Term] {\includegraphics[width=0.24\columnwidth]{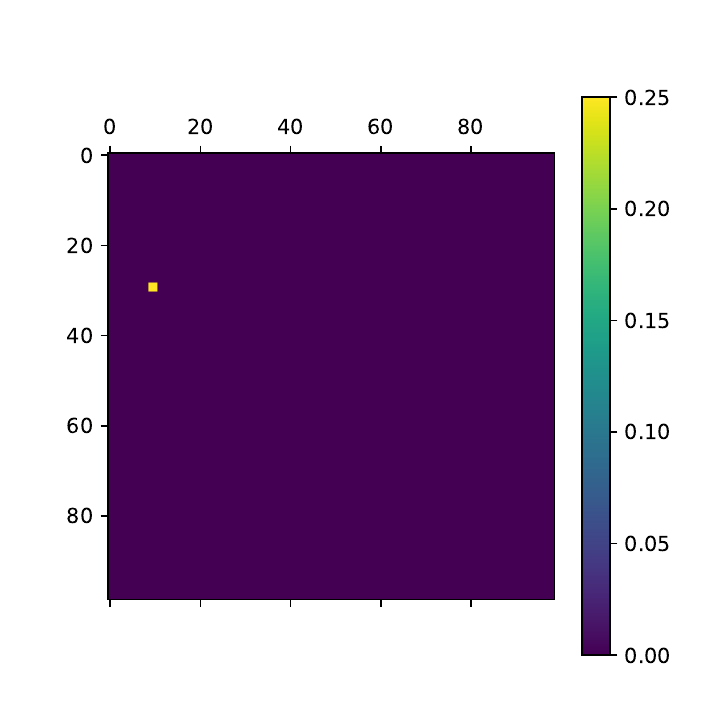}}
        \subfigure[Homogenized $\kappa^\star_{11}$] {\includegraphics[width=0.24\columnwidth]{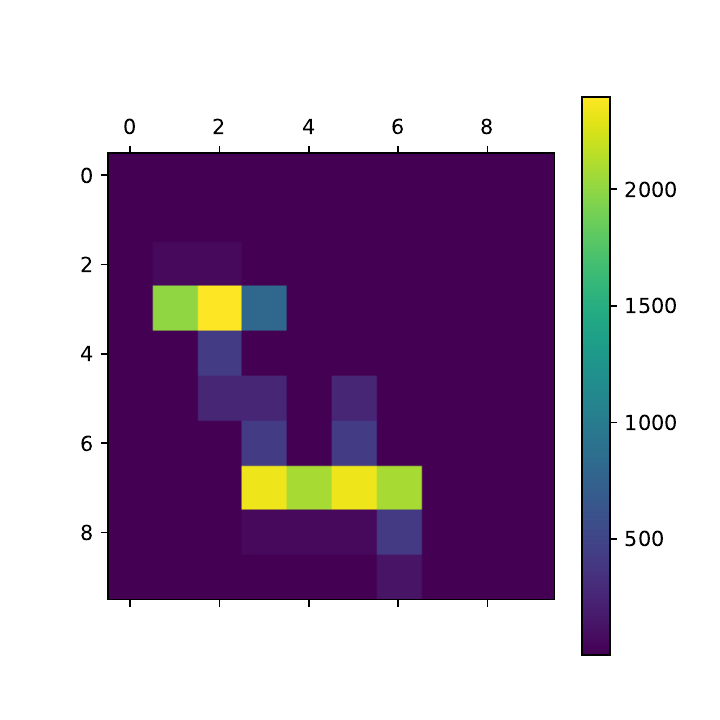}}
        \subfigure[Homogenized $\kappa^\star_{22}$] {\includegraphics[width=0.24\columnwidth]{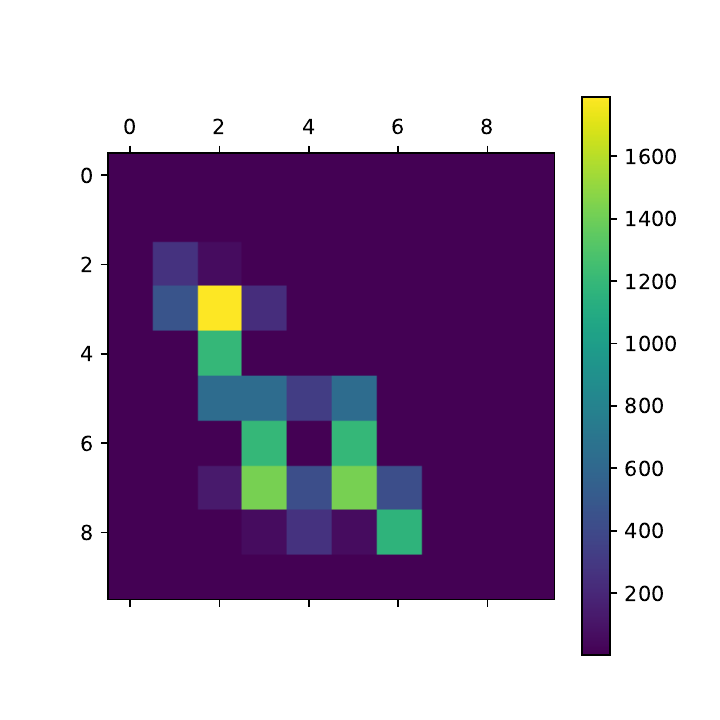}}
        \caption{The configuration of the permeability field and the point source term. (a) The permeability field. (b) The point source term. (c) The homogenized permeability $\kappa^\star_{11}$. (d) The homogenized permeability $\kappa^\star_{22}$.}
        \label{nonlinear-permeability}
\end{figure}

In the example of the nonlinear equation, the permeability field $\kappa(\mathbf{x},u(\mathbf{x},t))=a(\mathbf{x})e^{\beta u(\mathbf{x},t)}$, and $\beta$ is set to 0.01, which has relatively weak nonlinearity. The source term $f(\mathbf{x},t)$ is a point source term. The configuration of the permeability field and the source term can be found in \textbf{Figure} \ref{nonlinear-permeability}. The conductivity is $10^4$ in the channel and 1 in the background. Similar to before, we show the diagonal elements $\kappa^\star_{11}$ and $\kappa^\star_{22}$ of the homogenized permeability in \textbf{Figure} \ref{nonlinear-permeability} (c-d). Note we only show the permeability at time $t=0$, because the equation is relatively weakly nonlinear and the permeability changes with time are very small. The total simulation time $T=0.01$ and the time step size $\tau=0.001$.

\begin{figure}[t]
        \centering
        \subfigure[Reference solution $t=0.001$]
        {\includegraphics[width=0.25\columnwidth]{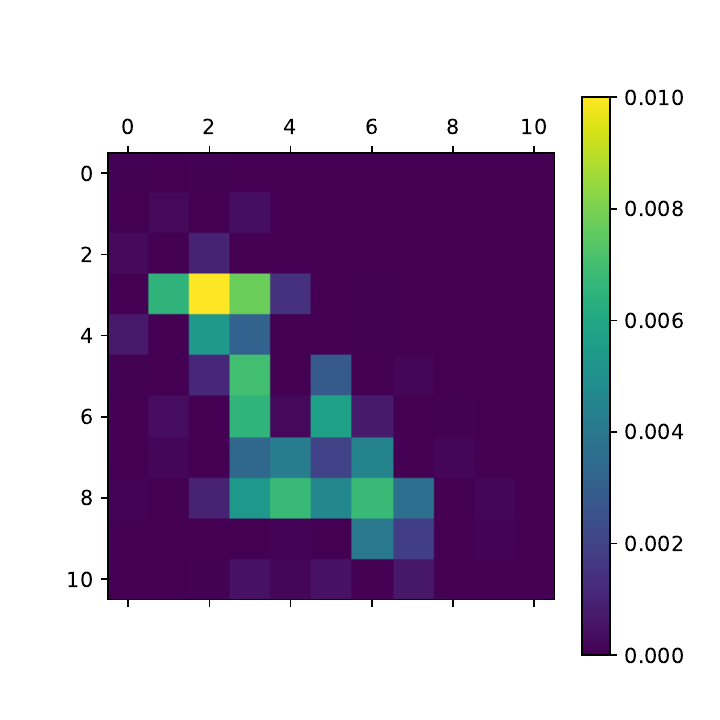}}
        \subfigure[Homogenized solution $t=0.001$] {\includegraphics[width=0.25\columnwidth]{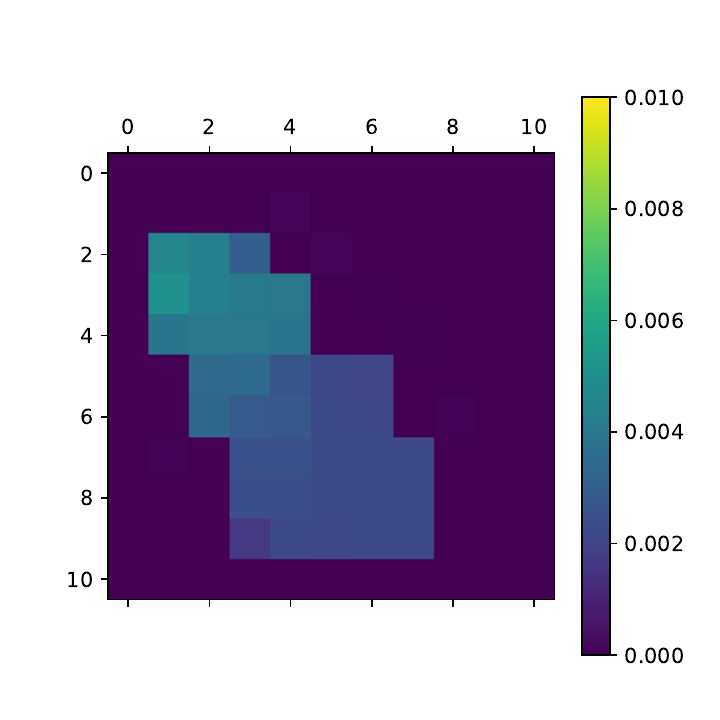}}
        \subfigure[Corrected solution $t=0.001$]
        {\includegraphics[width=0.25\columnwidth]{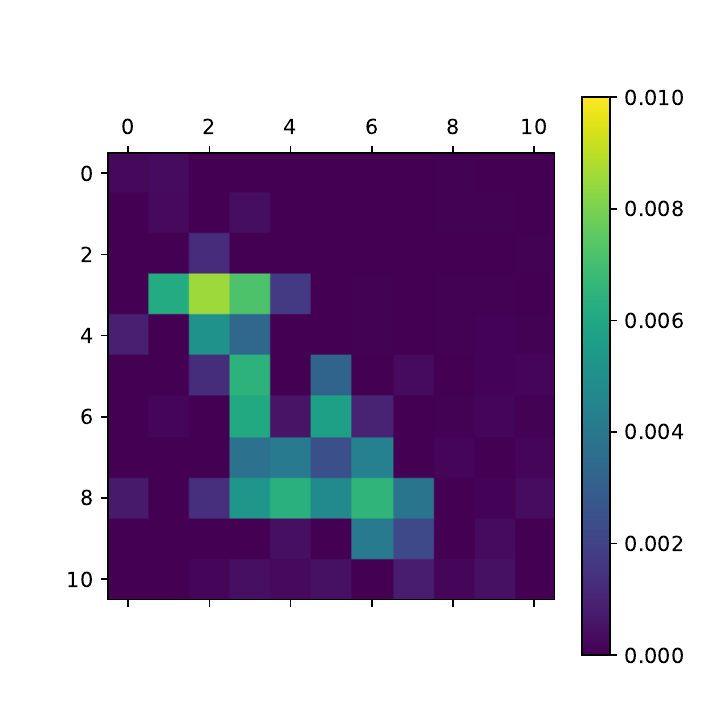}}\\

        \subfigure[Reference solution $t=0.01$]
        {\includegraphics[width=0.25\columnwidth]{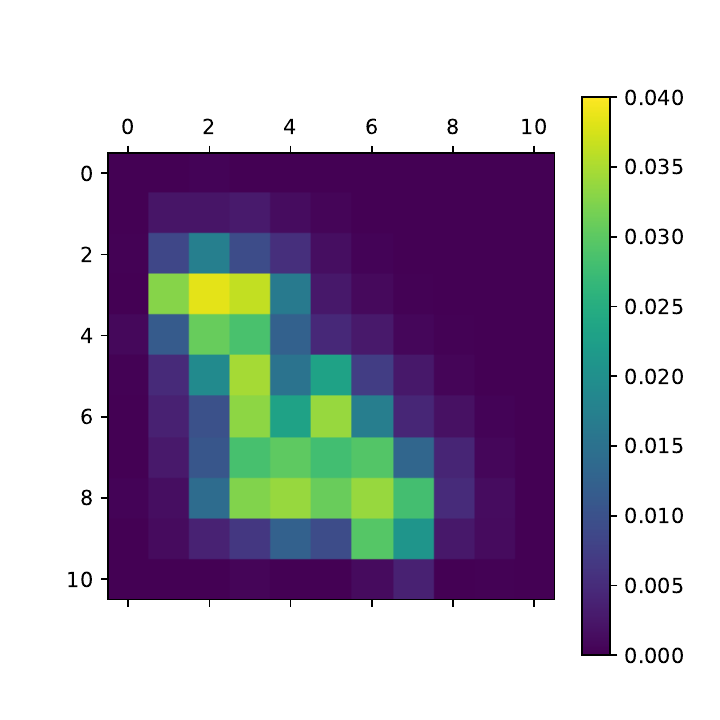}}
        \subfigure[Homogenized solution $t=0.01$] {\includegraphics[width=0.25\columnwidth]{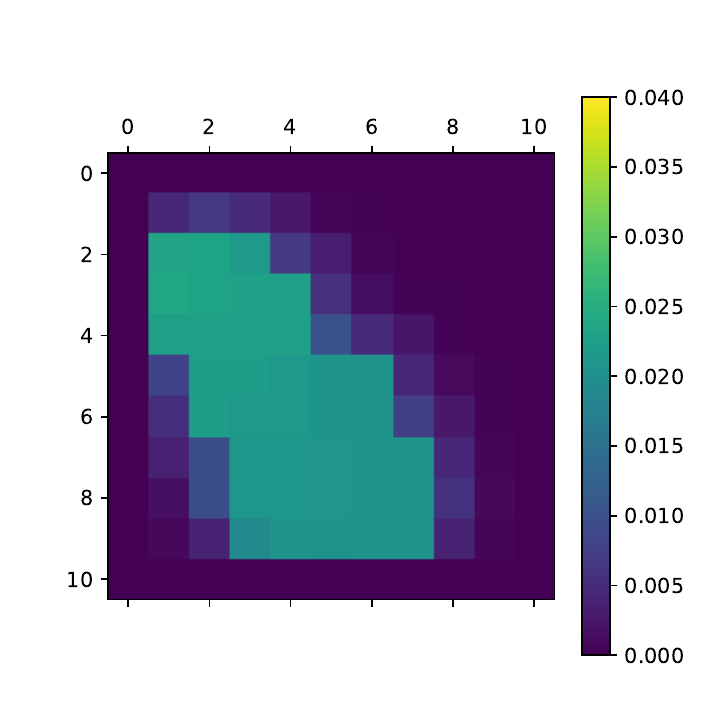}}
        \subfigure[Corrected solution $t=0.01$]
        {\includegraphics[width=0.25\columnwidth]{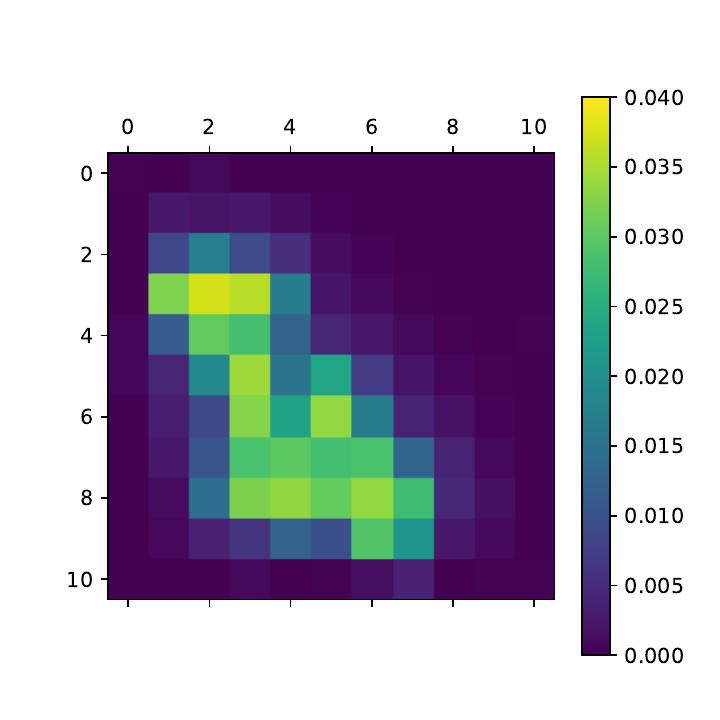}}

        \caption{The FEM solutions of the nonlinear flow equation, the homogenized equation and our proposed learning-based multi-continuum equation. (a)(d) The reference solutions at $t=0.001,0.01$; (b)(e) The homogenized solutions at $t=0.001,0.01$; (c)(f) The corrected solutions at $t=0.001,0.01$;.}
        \label{nonlinear-solution}
\end{figure}

The PINN is used as the PDE solver for the nonlinear flow equation, and network parameters are optimized through direct BP. Networks are trained for 300000 epochs and the loss function (\ref{pinn}) can reach $O(10^{-3})$. \textbf{Figure} \ref{nonlinear-solution} shows examples of the reference solutions, homogenized solutions, and the corrected solutions of our multi-continuum model. It can be seen that numerical homogenization can capture the overall dynamics of flow, with relatively large errors. In contrast, our learning-based multi-continuum model corrects the homogenized equation and results in more accurate solutions. \textbf{Table} \ref{tbl3} lists the relative $L_2$ error of our corrected solutions and homogenized solutions with respect to reference solutions at each time step. The corrected solutions obtained by our multi-continuum model are very close to the true solution, with the relative $L_2$ errors all within $1\%$. While the relative $L_2$ error of the numerical homogenization reaches $57.4263\%$ at time 0.001.

\begin{figure}[h]
        \centering
        \subfigure[Learned $\kappa_2^{11}(t=0.001)$]
        {\includegraphics[width=0.24\columnwidth]{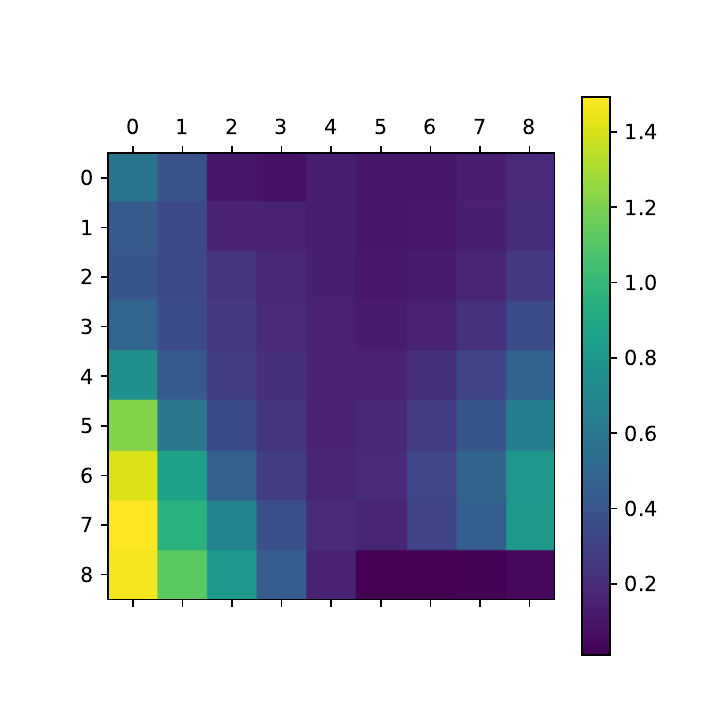}}
        \subfigure[Learned $\kappa_2^{11}(t=0.01)$]
        {\includegraphics[width=0.24\columnwidth]{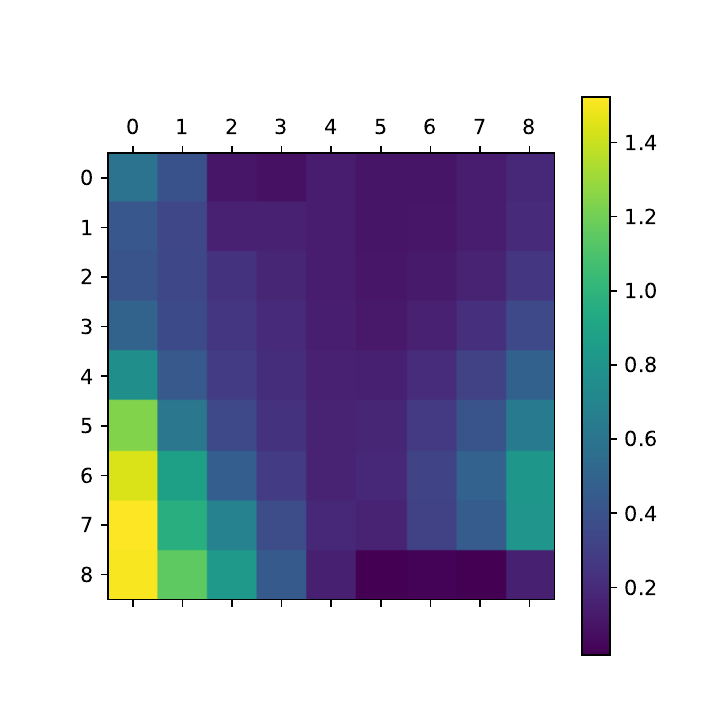}}
        \subfigure[Learned $\kappa_2^{22}(t=0.001)$]
        {\includegraphics[width=0.24\columnwidth]{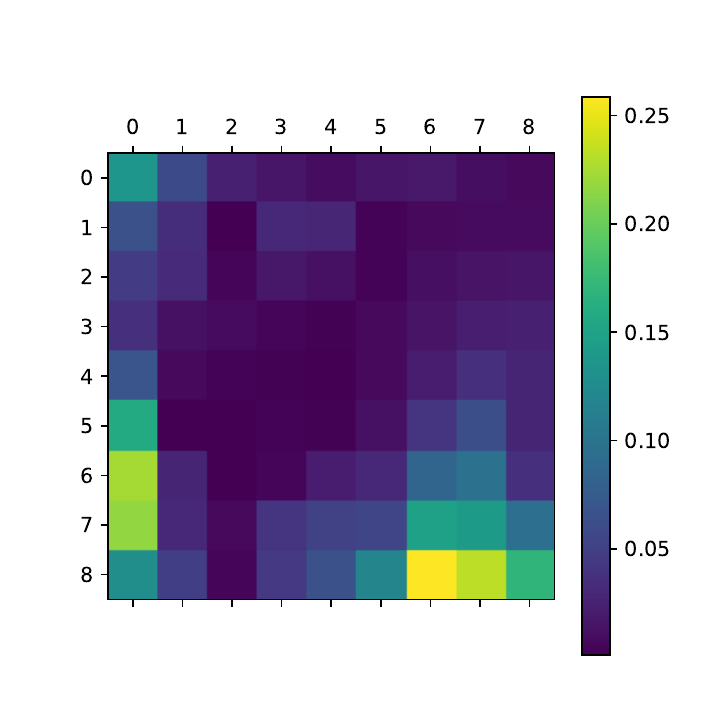}}
        \subfigure[Learned $\kappa_2^{22}(t=0.01)$]
        {\includegraphics[width=0.24\columnwidth]{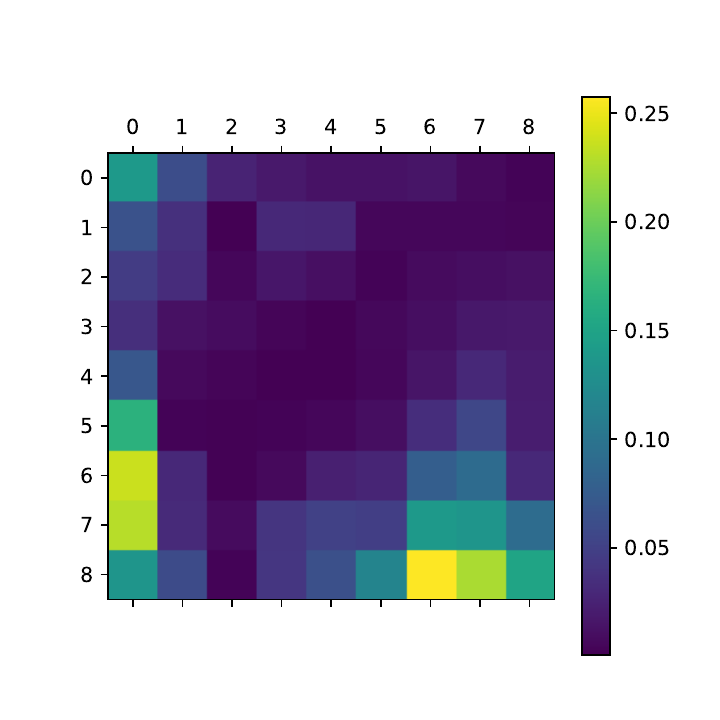}}
        \subfigure[Learned $\sigma(t=0.001)$]
        {\includegraphics[width=0.24\columnwidth]{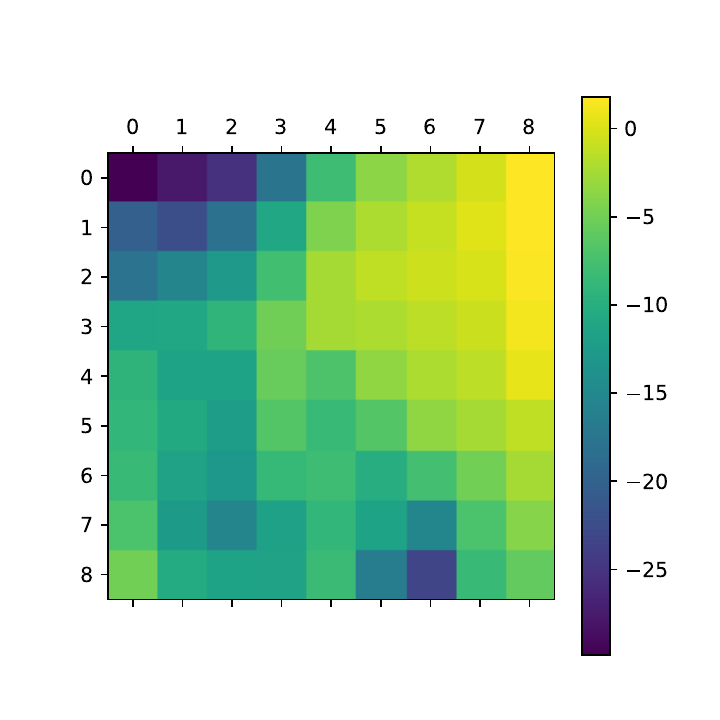}}
        \subfigure[Learned $\sigma(t=0.01)$]
        {\includegraphics[width=0.24\columnwidth]{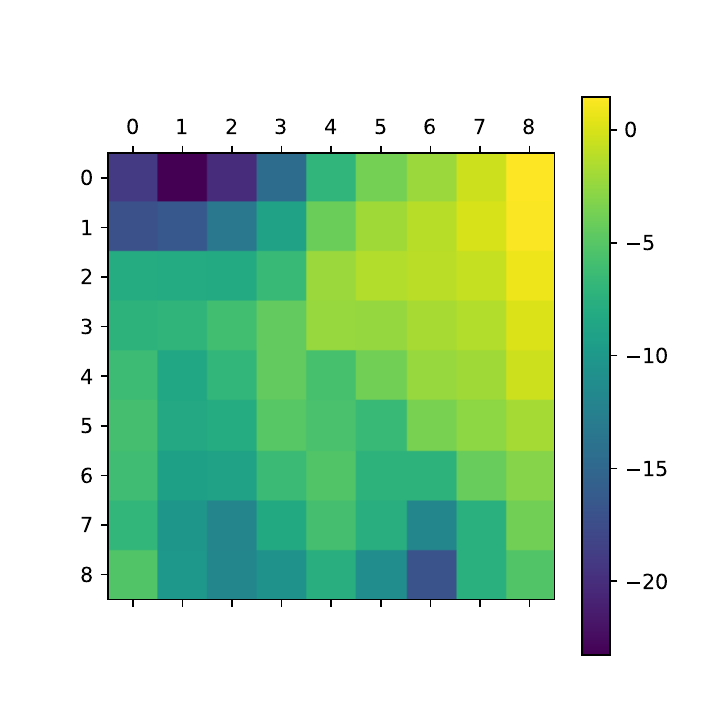}}
        \caption{The learned permeability field $\kappa_2$ and transfer coefficient $\sigma$. (a-d) The diagonal elements $\kappa_2^{11}$ and $\kappa_2^{22}$ at time $t=0.001$ and 0.01 in domain $\Omega$. (e-f) The learned transfer coefficient $\sigma$ at time $t=0.001$ and 0.01.}
        \label{nonlinear-learned-parameter}
\end{figure}

\textbf{Figure} \ref{nonlinear-learned-parameter} further illustrates the values of the learned permeability $\kappa_2(\mathbf{x},t;\theta_1)$ and the learned transfer coefficient $\sigma(\mathbf{x},t;\theta_2)$ for the nonlinear flow equation, which are all dependent on the time $t$. The diagonal elements $\kappa_2^{11}$ and $\kappa_2^{22}$ in the whole domain $\Omega$ at time $t=0.001$ and 0.01 are shown in \textbf{Figure} \ref{nonlinear-learned-parameter} (a-d). It can be observed that $\kappa_2^{11}$ and $\kappa_2^{22}$ vary very little over time and the values of $\kappa_2^{11}$ and $\kappa_2^{22}$ are also closely grouped across the domain $\Omega$. These results further indicate that the newly introduced continuum in the multi-continuum model can complement the information obtained from the numerical homogenization. The learned coefficient $\sigma$ at time 0.001 and 0.01 is illustrated in \textbf{Figure} \ref{nonlinear-learned-parameter} (e-f). It can be seen that $\sigma$ has essentially the same trend across the region at different time steps. The absolute values of $\sigma$ are larger all in places where the flow fluctuates more.

\begin{table}[p]
\caption{The relative $L_2$ error (in percentage) of the homogenized solutions and solutions of the proposed multi-continuum model with respect to the reference solutions.}\label{tbl3}
\begin{tabular*}{\tblwidth}{@{}LLLLLLLLLLL@{}}
\toprule
 \diagbox{Method}{Time} & 0.001 & 0.002 & 0.003 & 0.004 & 0.005 & 0.006 & 0.007 & 0.008 & 0.009 & 0.01 \\ % Table header row
\midrule
 Homogenization &  57.4263 & 50.9699 & 46.1707 & 42.4819 & 39.5630 & 37.1998 & 35.2499 & 33.6150 & 32.2250 & 31.0292\\
 Multi-continuum model &  0.6265 & 0.2800 & 0.0993 & 0.0343 & 0.0169 & 0.0132 & 0.0105 & 0.0061  & 0.0037 & 0.0124\\
\bottomrule
\end{tabular*}
\end{table}

\section{Conclusion}\label{conclusion}
In this paper, we propose a novel learning-based multi-continuum model for complex multiscale problems to correct solutions obtained through numerical homogenization and improve their accuracy. Without loss of generalization, we consider a two-continuum case. We reformulate the homogenized equation as a multi-continuum model, with the flow equation for the first continuum retaining the information of the original homogenized equation, along with an additional interaction term. The second continuum's flow equation necessitates the use of a neural network to determine its effective permeability. The interaction term between the two continua is aligned with that of the Dual-porosity model but with a learnable coefficient determined by another neural network. Using trusted data, i.e., accurate solutions at some points, we optimize these two neural networks by minimizing the mismatch between the numerical solutions of the multi-continuum model and the trusted data. We further discuss both direct back-propagation and the adjoint method for the PDE-constraint optimization problem. Our proposed learning-based multi-continuum model can resolve multiple interacted media within each coarse grid block and describe the mass transfer among them by introducing a new continuum to complement the missing/over-estimated information of the numerical homogenization. Numerical experiments demonstrate that our proposed learning-based multi-continuum model can correct the solutions obtained through numerical homogenization and improve their accuracy.

\section*{Acknowledgments}
This work was supported by National Natural Science Foundation of China (grant no. 12301559) and Hong Kong RGC Early Career Scheme (grant no. 9048274).

% \begin{figure}[t]
%         \centering
%         \subfigure[Reference solution $t=0.001$]
%         {\includegraphics[width=0.25\columnwidth]{nonlinear2_t2_sol.eps}}
%         \subfigure[Homogenized solution $t=0.001$] {\includegraphics[width=0.25\columnwidth]{nonlinear2_t6_sol.eps}}
%         \subfigure[Corrected solution $t=0.001$]
%         {\includegraphics[width=0.25\columnwidth]{nonlinear2_t11_sol.eps}}\\
%         {\includegraphics[width=0.25\columnwidth]{figures/nonlinear2_errs.eps}}\\
%         \caption{The FEM solutions of the nonlinear flow equation, the homogenized equation and our proposed learning-based multi-continuum equation. (a)(d) The reference solutions at $t=0.001,0.01$; (b)(e) The homogenized solutions at $t=0.001,0.01$; (c)(f) The corrected solutions at $t=0.001,0.01$;.}
%         \label{nonlinear2-solution}
% \end{figure}


\section*{Appendix}
\appendix
\textbf{A. Derivation of adjoint equation.}
The loss function is given by
\begin{equation}\label{continuous}
 \mathcal{L} = \frac{1}{2}\int_0^T \int_{\Omega} \| u_1(\mathbf{x},t) - U(\mathbf{x},t) \|^2dxdt.
\end{equation}

The first step is to introduce the Lagrangian corresponding to the optimization problem:

\begin{equation}
\begin{aligned}
 \mathcal{L} \equiv & \int_0^T \int_{\Omega} \frac{1}{2}\| u_1(\mathbf{x},t) - U(\mathbf{x},t) \|^2 + \lambda_1(\mathbf{x},t) \big[\frac{\partial u_1(\mathbf{x}, t)}{\partial t} - \text{div}(\kappa^\star_1(\mathbf{x})\nabla u_1(\mathbf{x},t)) + \sigma(\mathbf{x};\theta_2)(u_1(\mathbf{x},t)-u_2(\mathbf{x},t))-f(\mathbf{x},t)\big] \\
  & + \lambda_2(\mathbf{x},t) \big[\frac{\partial u_2(\mathbf{x,t})}{\partial t} - \text{div}(\kappa_2(\mathbf{x};\theta_1)\nabla u_2(\mathbf{x},t)) + \sigma(\mathbf{x};\theta_2)(u_2(\mathbf{x},t)-u_1(\mathbf{x},t))\big]dxdt + \int_{\Omega} \mu_1 u_1(\mathbf{x},0) + \mu_2 u_2(\mathbf{x},0)dx \\
  & + \int_0^T \int_{\partial \Omega} \mu_3 u_1(\mathbf{x},t) + \mu_4  u_2(\mathbf{x},t) dxdt.
\end{aligned}
\end{equation}
The Lagrangian multipliers $\lambda_1$ and $\lambda_2$ is a function of time $t$ and coordinate $x$, and $\mu_1, \mu_2,\mu_3,\mu_4$ are other multipliers that are associated with the initial and boundary conditions.

\begin{equation}\label{continuous-grdient}
\begin{aligned}
  \frac{\partial \mathcal{L}}{\partial \theta_1} = &  \int_0^T \int_{\Omega} (u_1(\mathbf{x},t) - U(\mathbf{x},t) ) \frac{\partial u_1}{\partial \theta_1}dxdt \\
  & + \int_0^T \int_{\Omega}\lambda_1 \big[\frac{\partial^2 u_1(\mathbf{x}, t)}{\partial \theta_1 \partial t} -  \frac{\partial}{\partial \theta_1} \text{div}(\kappa^\star_1(\mathbf{x})\nabla u_1(\mathbf{x},t)) + \sigma(\mathbf{x};\theta_2)(
 \frac{u_1(\mathbf{x},t)}{\partial \theta_1}-\frac{u_2(\mathbf{x},t)}{\partial \theta_1})\big]dxdt \\
  & + \int_0^T \int_{\Omega}\lambda_2 \big[\frac{\partial^2 u_2(\mathbf{x},t)}{\partial \theta_1\partial t} - \frac{\partial}{\partial \theta_1} \text{div}(\kappa_2(\mathbf{x};\theta_1)\nabla u_2(\mathbf{x},t)) + \sigma(\mathbf{x};\theta_2)(\frac{u_2(\mathbf{x},t)}{\partial \theta_1}-\frac{u_1(\mathbf{x},t)}{\partial \theta_1})\big]dxdt \\
  & + \int_{\Omega} \mu_1 \frac{\partial u_1(\mathbf{x},0)}{\partial \theta_1}+ \mu_2 \frac{\partial u_2(\mathbf{x},0) }{\partial \theta_1}dx + \int_0^T \int_{\partial \Omega} \mu_3 \frac{\partial u_1(\mathbf{x},t)}{\partial \theta_1} + \mu_4  \frac{\partial u_2(\mathbf{x},t) }{\partial \theta_1} dxdt.
\end{aligned}
\end{equation}
The next step is to integrate by parts to eliminate the terms $\frac{\partial^2 u_1(\mathbf{x}, t)}{\partial \theta_1 \partial t}$ and $\frac{\partial^2 u_2(\mathbf{x}, t)}{\partial \theta_1 \partial t}$:
\begin{equation}
\begin{aligned}
  &\int_0^T \int_{\Omega} \lambda_1 \big[\frac{\partial^2 u_1(\mathbf{x}, t)}{\partial \theta_1 \partial t} -  \frac{\partial}{\partial \theta_1} \text{div}(\kappa^\star_1(\mathbf{x})\nabla u_1(\mathbf{x},t))\big]dxdt \\
  &=\int_{\Omega} (\lambda_1\frac{\partial u_1(\mathbf{x}, t)}{\partial\theta_1}\big|_0^T - \int_0^T \frac{\partial u_1(\mathbf{x}, t)}{\partial\theta_1} \frac{\partial \lambda_1}{\partial t}dt)dx - \int_0^T \int_{\Omega} \frac{\partial u_1(\mathbf{x}, t)}{\partial\theta_1} \text{div}(\kappa^\star_1(\mathbf{x})\nabla \lambda_1)dxdt \\
  &=\int_0^T \int_{\Omega} -\frac{\partial u_1(\mathbf{x}, t)}{\partial\theta_1} \frac{\partial \lambda_1}{\partial t} - \frac{\partial u_1(\mathbf{x}, t)}{\partial\theta_1} \text{div}(\kappa^\star_1(\mathbf{x})\nabla \lambda_1) dxdt - \int_\Omega \lambda_1 (\mathbf{x},0) \frac{\partial u_1(\mathbf{x}, 0)}{\partial\theta_1}dx,
\end{aligned}
\end{equation}
where the transformation in the second and third line use the Divergence Theorem, and we let the boundary condition $\lambda_1(\mathbf{x},t)|_{\partial \Omega} = 0$ and final condition $\lambda_1(\mathbf{x}, T)=0$. The same result can be obtained for $u_2$, i.e.,
\begin{equation}
\begin{aligned}
  & \int_0^T \int_{\Omega}  \lambda_2 \big[\frac{\partial^2 u_2(\mathbf{x,t})}{\partial \theta_1\partial t} - \frac{\partial}{\partial \theta_1} \text{div}(\kappa_2(\mathbf{x};\theta_1)\nabla u_2(\mathbf{x},t))\big]dxdt \\
  = & \int_0^T \int_{\Omega} -\frac{\partial u_2(\mathbf{x}, t)}{\partial\theta_1} \frac{\partial \lambda_2}{\partial t} - \frac{\partial u_2(\mathbf{x}, t)}{\partial\theta_1} \text{div}(\kappa_2(\mathbf{x}, \theta_1)\nabla \lambda_2) dxdt\\
  & - \int_0^T \int_{\Omega} \lambda_2 \text{div}( \frac{\kappa_2(\mathbf{x}, \theta_1)}{\partial \theta_1}\nabla u_2(\mathbf{x},t))dxdt - \int_\Omega \lambda_2 (\mathbf{x},0) \frac{\partial u_2(\mathbf{x}, 0)}{\partial\theta_1}dx.
  \end{aligned}
\end{equation}

Substituting these results into equation (\ref{continuous-grdient}),

\begin{equation}
\begin{aligned}
\frac{\partial \mathcal{L}}{\partial \theta_1} = & \int_0^T \int_{\Omega} ( u_1(\mathbf{x},t) - U(\mathbf{x},t) ) \frac{\partial u_1}{\partial \theta_1}dxdt \\
 & + \int_0^T \int_{\Omega} -\frac{\partial u_1(\mathbf{x}, t)}{\partial\theta_1} \frac{\partial \lambda_1}{\partial t} - \frac{\partial u_1(\mathbf{x}, t)}{\partial\theta_1} \text{div}(\kappa^\star_1(\mathbf{x})\nabla \lambda_1) +\lambda_1 \sigma(\mathbf{x};\theta_2)( \frac{u_1(\mathbf{x},t)}{\partial \theta_1}-\frac{u_2(\mathbf{x},t)}{\partial \theta_1})dxdt \\
 & + \int_0^T \int_{\Omega} -\frac{\partial u_2(\mathbf{x}, t)}{\partial\theta_1} \frac{\partial \lambda_2}{\partial t} - \frac{\partial u_2(\mathbf{x}, t)}{\partial\theta_1} \text{div}(\kappa_2(\mathbf{x}, \theta_1)\nabla \lambda_2) - \lambda_2 \text{div}( \frac{\kappa_2(\mathbf{x}, \theta_1)}{\partial \theta_1}\nabla u_2(\mathbf{x},t))dxdt\\
 & + \int_0^T \int_{\Omega} \lambda_2 \sigma(\mathbf{x};\theta_2)(\frac{u_2(\mathbf{x},t)}{\partial \theta_1}-\frac{u_1(\mathbf{x},t)}{\partial \theta_1})dxdt
  + \int_{\Omega} \mu_1 \frac{\partial u_1(\mathbf{x},0)}{\partial \theta_1}+ \mu_2 \frac{\partial u_2(\mathbf{x},0) }{\partial \theta_1}dx \\
  & - \int_\Omega \lambda_1 (\mathbf{x},0) \frac{\partial u_1(\mathbf{x}, 0)}{\partial\theta_1}dx - \int_\Omega \lambda_2 (\mathbf{x},0) \frac{\partial u_2(\mathbf{x}, 0)}{\partial\theta_1}dx + \int_0^T \int_{\partial \Omega} \mu_3 \frac{\partial u_1(\mathbf{x},t)}{\partial \theta_1}+ \mu_4 \frac{\partial u_2(\mathbf{x},t) }{\partial \theta_1} dxdt.
\end{aligned}
\end{equation}
We set $\mu_1 = \lambda_1 (\mathbf{x},0)$ and $\mu_2 = \lambda_2 (\mathbf{x},0) $, $\mu_3=\mu_4=0$ to cancel the last four terms. The calculation of gradient $\frac{\partial u_1(\mathbf{x}, t)}{\partial\theta_1}$ and $\frac{\partial u_2(\mathbf{x}, t)}{\partial\theta_1}$ can be avoided by setting

\begin{equation}\label{adjoint-equation-appedx}
\begin{aligned}
 & \frac{\partial \lambda_1}{\partial t} + \text{div}(\kappa^\star_1(\mathbf{x})\nabla \lambda_1) - \sigma(\mathbf{x}; \theta_2)(\lambda_1 - \lambda_2) =  u_1(\mathbf{x},t) - U(\mathbf{x},t),  \\
 & \frac{\partial \lambda_2}{\partial t} + \text{div}(\kappa_2(\mathbf{x}, \theta_1)\nabla \lambda_2) - \sigma(\mathbf{x}; \theta_2)(\lambda_2 - \lambda_1)= 0.
 \end{aligned}
\end{equation}
This is the obtained adjoint equation with the initial condition $\lambda_1(\mathbf{x},t=T)=\lambda_2(\mathbf{x},t=T)=0$ and boundary condition $\lambda_1(\mathbf{x},t)|_{\partial \Omega}=\lambda_2(\mathbf{x},t)|_{\partial \Omega}=0$. Finally the gradient can be obtained by
\begin{equation}
\begin{aligned}
 \frac{\partial \mathcal{L}}{\partial \theta_1}  = \int_0^T \int_{\Omega} -\lambda_2 \text{div}( \frac{\kappa_2(\mathbf{x}, \theta_1)}{\partial \theta_1}\nabla u_2(\mathbf{x},t))dxdt.
\end{aligned}
\end{equation}
The similar derivation can be applied for $\theta_2$ and the obtained adjoint equation is the same as equation (\ref{adjoint-equation-appedx}). The desired gradient are given by

\begin{equation}
\begin{aligned}
 \frac{\partial \mathcal{L}}{\partial \theta_2}  = \int_0^T \int_{\Omega} \lambda_1 \frac{\partial \sigma(\mathbf{x};\theta_2)}{\partial \theta_2} (u_1(\mathbf{x},t)-u_2(\mathbf{x},t)) + \lambda_2 \frac{\partial \sigma(\mathbf{x};\theta_2)}{\partial \theta_2} (u_2(\mathbf{x},t)-u_1(\mathbf{x},t))dxdt.
\end{aligned}
\end{equation}
In this way, the network parameters $\theta_1$ and $\theta_2$ can be optimized by classical optimization method, such as gradient descent method. In the integral formulations, (.8) and (.9), the derivatives  $\frac{\kappa_2(\mathbf{x}, \theta_1)}{\partial \theta_1}$ and $\frac{\partial \sigma(\mathbf{x};\theta_2)}{\partial \theta_2}$ can be calculated by the automatic differentiation of Pytorch/Tensorflow.

Further, we calculate the numerical solutions of adjoint equation by FEM. The corresponding discretization of adjoint equation is given by
\begin{equation}
\begin{aligned}
 & M\frac{\partial \lambda_1}{\partial t} - A_1(\kappa_1^\star)\lambda_1 -\sigma(\theta_2)M(\lambda_1-\lambda_2) = M(u-U),\\
 & M\frac{\partial \lambda_2}{\partial t} - A_2(\kappa_2(\theta_1))\lambda_2 -\sigma(\theta_2)M(\lambda_2-\lambda_1) = 0,
\end{aligned}
\end{equation}
where $A_1,A_2$ are the stiffness matrices matrices, and $M$ is the mass matrix. The gradient of the loss function is reformulated as
\begin{equation}
\begin{aligned}
 \nabla_{\theta_1,\theta_2} \mathcal{L} =&  \int_T^0  \lambda_2^T \nabla_{\theta_1} A_2(\kappa_2)u_2dt\\
 +& \int_T^0 \lambda_1^T \nabla_{\theta_2} \sigma(\theta_2)M(u_1-u_2)dt\\
 +& \int_T^0 \lambda_2^T \nabla_{\theta_2} \sigma(\theta_2)M(u_2-u_1)dt.
\end{aligned}
\end{equation}

\txrev{\textbf{B. The details on the trusted data and the reference solutions}. The reference solutions are obtained by averaging the fine-scale solutions over each coarse block. In the above experiments in Section~\ref{experiments}, the trusted data used to train neural networks are the whole reference solutions. This is because our proposed method is designed not only to correct the homogenized solution but also to obtain a novel multi-continuum model, which more accurately characterizes the intricate dynamics within fractured media. Here we conduct further experiments to evaluate the performance of our proposed method when the trusted data is sampled from the reference solutions.}

\txrev{The initial stage of the experiments involves sampling the trusted data in time from the reference solution. The results are presented in Table \ref{tbl4}. It can be observed that the relative $L_2$ error exhibits a slight reduction with more time steps sampled (from $t=0.4$ to $t=0.6$). This outcome reflects that as the number of sampled time steps increases, the model is able to characterize the dynamics with greater accuracy and obtain more precise solutions. Furthermore, the observed reduction in error is relatively modest, indicating that our method is capable of learning a robust model from a limited amount of time-step data.}

\txrev{Furthermore, the trusted data were randomly sampled from the reference solutions at 60\%, 70\%, and 80\% of the total space, respectively, and all time steps were used. The results are listed in Table \ref{tbl5}. It can be seen that there is a significant increase in the relative $L_2$ error of the solution obtained by our model with fewer trusted data sampled in space. This indicates that our approach is reliant on the presence of a greater density of data in space. Furthermore, it can be postulated that the sampling strategy employed may also influence the results. In the current experiment, the trusted data were sampled from a uniform distribution. Should the sampling points be concentrated in the more variable regions of the dynamical system, it may be possible to achieve a performance improvement with less data. The impact of alternative sampling techniques, such as importance sampling, can be evaluated in future work.}

       \begin{table}[htp]
        \caption{The relative $L_2$ error (in percentage)  of the homogenized solutions and solutions of the proposed multi-continuum model with respect to the reference solutions when the trusted data are sampled in time.}\label{tbl4}
        \begin{tabular*}{\textwidth}{ccccccccccc}
        \toprule
         \diagbox{Method}{Time} & 0.1 & 0.2 & 0.3 & \multicolumn{1}{c}{0.4} & 0.5 & 0.6 & 0.7 & 0.8 & 0.9 & 1.0 \\
        \midrule
         Homogenization & 5.7956 & 7.0050 & 7.5187 & \multicolumn{1}{c}{7.7016} & 7.7560 & 7.7774 & 7.7824 & 7.7838 & 7.7841 & 7.7845\\
         \midrule
        \multirow{2}*{Our model} &  \multicolumn{4}{c|}{train} & \multicolumn{6}{c}{test} \\ \cline{2-11}
        &  0.1290 & 0.1285 & 0.1119 & \multicolumn{1}{c|}{0.1138} & 0.1145 & 0.1148 & 0.1149 & 0.1149 & 0.1149 & 0.1149\\
        \midrule
        \multirow{2}*{Our model}& \multicolumn{6}{c|}{train} & \multicolumn{4}{c}{test} \\ \cline{2-11}
         &  0.1251 & 0.1218 & 0.1130 & 0.1084 & 0.1137 & \multicolumn{1}{c|}{0.1142} & 0.1143 & 0.1143 & 0.1143 & 0.1143\\
        \bottomrule
        \end{tabular*}
        \end{table}

        \begin{table}[t]
        \caption{The relative $L_2$ error (in percentage) of the homogenized solutions and solutions of the proposed multi-continuum model with respect to the reference solutions when the trusted data are sampled in space. Sampling ratios are listed in parentheses.}\label{tbl5}
        \begin{tabular*}{\textwidth}{ccccccccccc}
        \toprule
         \diagbox{Method}{Time} &
          0.1 & 0.2 & 0.3 & 0.4 & 0.5 & 0.6 & 0.7 & 0.8 & 0.9 & 1.0 \\
        \midrule
         Homogenization & 5.7956 & 7.0050 & 7.5187 & 7.7016 & 7.7560 & 7.7774 & 7.7824 & 7.7838 & 7.7841 & 7.7845\\
        Our model (60\%)&  0.2415 & 0.2331 & 0.2303 & 0.2293 & 0.2290 & 0.2289 & 0.2288 & 0.2288 & 0.2288 & 0.2288\\
        Our model (70\%)&  0.1634 & 0.1630 & 0.1655 & 0.1669 & 0.1674 & 0.1676 & 0.1677 & 0.1677 & 0.1677 & 0.1677\\
        Our model (80\%)&  0.1533 & 0.1523 & 0.1554 & 0.1573 & 0.1581 & 0.1583 & 0.1584 & 0.1584 & 0.1584 & 0.1584\\
        Our model (100\%)&  0.1162 & 0.1148 & 0.1106 & 0.1091 & 0.1086 & 0.1085 & 0.1084 & 0.1084 & 0.1084 & 0.1084\\
        \bottomrule
        \end{tabular*}
        \end{table}

\printcredits

\bibliographystyle{cas-model2-names}
\bibliography{cas-refs.bib}

%% Biography
%% Here goes the biography details.
%% Here goes the biography details.

\end{document}